\documentclass[moor]{informs1}              % for a regular run
\usepackage{natbib}
 \NatBibNumeric
 \def\bibfont{\small}%
 \bibpunct[, ]{[}{]}{,}{n}{}{,}%

\usepackage[colorlinks=true,breaklinks=true,bookmarks=true,urlcolor=blue,
     citecolor=blue,linkcolor=blue,bookmarksopen=false,draft=false]{hyperref}

\def\EMAIL#1{\href{mailto:#1}{#1}}% When hyperref is used, otherwise outcomment 
\TheoremsNumberedThrough     % Preferred (Theorem 1, Lemma 1, Theorem 2)
\EquationsNumberedThrough    % Default: (1), (2), ...
\theoremstyle{TH}
\newtheorem{condition}{Condition}
\theoremstyle{TH}
\newtheorem{approximation}{Approximation}

\usepackage{multirow}
\usepackage[graphicx]{realboxes}
\usepackage{graphicx}
\usepackage{adjustbox}

\usepackage{comment}

\renewcommand{\theARTICLETOP}{}

\begin{document}
%%%%%%%%%%%%%%%%

%\setcounter{page}{1}
%\DOI{10.1287/xxxx.0000.0000}%

\RUNAUTHOR{Li and Mehrotra}

\RUNTITLE{A Finite Approximation Approach for Stochastic Model Analysis}

% \TITLE{Bundling Information Goods of Decreasing Value}

%\TITLE{A Near-Optimal Finite Approximation Approach for Computing Stationary Distributions and Performance Measures of Continuous-State Markov Chains}

%\TITLE{A Near Optimal Method for Computing Stationary Distributions and Performance Measures of Markov Chains using Finite Approximation}

%\TITLE{Computing Stationary Distributions and Performance Measures of Continuous-State Markov Chains using Finite Approximation with an Application to G/G/1 + G Queues}

%\TITLE{A Near Optimal Method for Computing Stationary Distributions and Performance Measures of Markov Chains using Finite Approximation with an Application to G/G/1+G Queues}

%\TITLE{A Near Optimal Method for Computing Stationary Distributions and Performance Measures of Markov Chains using Finite Approximation with a Queuing Application}

%\TITLE{A Near Optimal Method for Computing Stationary Distributions and Performance Measures of Markov Chains with a Queuing Application}

\TITLE{A New Method for Computing Stationary Distribution and Steady-State Performance Measures of a Continuous-State Markov Chain with a Queuing Application}

\ARTICLEAUTHORS{%
\AUTHOR{Shukai Li}
\AFF{Northwestern University, \EMAIL{shukaili2024@u.northwestern.edu}}
\AUTHOR{Sanjay Mehrotra}
\AFF{Northwestern University, \EMAIL{mehrotra@iems.northwestern.edu}}
} % end of the block

\ABSTRACT{%
Applications of stochastic models often involve the evaluation of steady-state performance, which requires solving a set of balance equations.
In most cases of interest, the number of equations is infinite or even uncountable. As a result, numerical or analytical solutions are unavailable.
This is true even when the system state is one-dimensional.
This paper develops a general method for computing stationary distributions and steady-state performance measures of stochastic systems that can be described as continuous-state Markov chains supported on $\mathbb{R}$. 
The balance equations are numerically solved by properly constructing a proxy Markov chain with finite states.  
We show the consistency of the approximate solution and provide deterministic non-asymptotic error bounds under the supremum norm. 
Our finite approximation method is near-optimal among all approximation methods using discrete distributions, including the empirical distributions generated by a simulation approach. We apply the developed method to compute the stationary distribution of virtual waiting time in a G/G/1+G queue and associated performance measures under certain mild but general differentiability and boundedness assumptions on the inter-arrival, service and patience time distributions. 
Numerical experiments validate the accuracy and efficiency of our method, and show it outperforms a standard Markov chain Monte Carlo method by several orders of magnitude. The developed method is also significantly more accurate than the available fluid approximations for this queue.
}%

%\KEYWORDS{deterministic inventory theory; infinite linear programming duality; 
%  existence of optimal policies; semi-Markov decision process; cyclic schedule}
%\MSCCLASS{Primary: 90B05; secondary: 90C40, 90C90}
%\ORMSCLASS{Primary: Inventory/production: deterministic multi-item;
%  secondary: dynamic programming/optimal control: deterministic 
%  semi-Markov; programming: infinite dimensional}
%\HISTORY{Received November 20, 2003; revised March 8, 2004, and March 26, 2004.}

\KEYWORDS{performance evaluation, stationary  distribution, finite approximation, stochastic modeling.}
\MSCCLASS{}
\ORMSCLASS{Primary: ; secondary: }
\HISTORY{}

\maketitle
%%%%%%%%%%%%%%%%%%%%%%%%%%%%%%%%%%%%%%%%%%%%%%%%%%%%%%%%%%%%%%%%%%%%%%

% Samples of sectioning (and labeling) in MOOR.
% NOTE: (1) all section levels end with a period,
%       (2) capitalization is as shown (sentence style, not title style).
%
%\section{Introduction.}\label{intro} %%1.
%\subsection{Duality and the classical EOQ problem.}\label{class-EOQ} %% 1.1.
%\subsection{Outline.}\label{outline1} %% 1.2.
%\subsubsection{Cyclic schedules for the general deterministic SMDP.}
%  \label{cyclic-schedules} %% 1.2.1
%\section{Problem description.}\label{problemdescription} %% 2.

% Text of your paper here
\section{Introduction.}
%\subsection{Overview of the Problem}
%Stochastic models that can be described as Markov chains have important applications.

Many stochastic systems can be described by a Markov chain (MC).
Computing steady-state performance measures for such systems involves solving a set of \textit{balance equations}. 
As pointed out by \cite{kuntz2021stationary}, ``in most cases of interest, the number of equations is infinite or too large and the equations cannot be solved analytically or numerically." 
This is true even when the system state is one-dimensional. 
Such systems have numerous applications in operations \cite{daley1965general, baccelli1981queues, zeltyn2005call, ward2008asymptotically, bassamboo2010accuracy, bremaud2013markov}, finance and economics \cite{tauchen1986finite,duan2001american,nicolau2002stationary, kopecky2010finite}.

This paper considers the steady-state performance evaluation problem of a (discrete time) MC with one-dimensional states.  
%Markov chain (MC) models often require knowledge of the stationary distribution and evaluation of system performance; e.g., . 
%In many cases of interest, the state space is continuous, mixed, or finite but too large, and analytical solutions to stationary  distribution equations  are unavailable; see \cite{kuntz2021stationary} for a review. In these situations, computation is often performed using the Monte Carlo approaches. 
Such a problem can be formulated as:
\begin{align}
    \mathbb{E}[\phi(Z)]=&\int_{\Omega} \phi(z) dp(z),\label{eqn:core problem 1}\tag{P1} \\ 
    \text{where}\quad\, p(x)\,\,\,= &\int_{\Omega}\tau(x,u)dp(u),\quad \forall x\in\Omega. \label{eqn:core problem 2}\tag{P2}
\end{align}
Here $\Omega\subseteq \mathbb{R}$ is a support, $z\in\Omega$ is the system state, $\phi:\Omega\mapsto\mathbb{R}$ is a known performance function, $\tau:\Omega^2\mapsto\mathbb{R}_+$ is defined by the MC's transition kernel, and $p$ is the \textit{stationary (probability) distribution} that can be obtained by solving the balance equations \eqref{eqn:core problem 2}; see an example in Subsection \ref{sec:GG1G}. 
We assume  $\Omega=[0,1]$ for the general results of this paper. However, our analysis can be adapted for other support $\Omega'\subseteq \mathbb{R}$  using certain transformations (see Section \ref{sec:prelim}).

We develop a deterministic finite approximation approach for \eqref{eqn:core problem 1}-\eqref{eqn:core problem 2}. Specifically, we approximate the original MC with properly constructed and increasingly refined finite-state MCs, whose stationary distributions are computable. 
The idea of finite approximation is natural, whereas
%In practice, however, the use of numerical integration is limited except for some seldom cases where stationary  distributions have analytical expressions. Indeed, stationary  distributions are often in an implicit form of integral equations that are  hard to solve particularly when the system states are continuous or mixed%. 
%A natural idea for overcoming the obstacle of solving the continuous-state stationary distribution, and performing numerical integration for the performance measure computations is to approximate 
the development of a general approximation scheme with tractable error bounds is not known. 
Existing methods are often restricted to certain types of models, %and how to generalize for other models is unclear, 
e.g., \cite{baumann2013computing, dayar2012kronecker,kazeev2015tensor, kuntz2019bounding}, while others focus on countable-state systems, e.g., \cite{zhao1996censored, tweedie1998truncation, masuyama2016error, masuyama2016limit, masuyama2017continuous, masuyama2017error, liu2018error,liu2018error2}. 
The accuracy of existing methods is often validated by experiments only or asymptotic results, e.g., \cite{pincus1992approximating, hart2012convergence, herve2019state}, but the computation of errors is not known, which limits their broader application. 
Thus, the development of a general finite approximation scheme with computable error bounds, as accomplished in this paper, fills an important fundamental gap in research on stochastic systems.

\subsection{Contributions and Techniques.}\label{sec:contributions}
\, \medskip \\
\noindent \textbf{A General and Computable Theory.} \quad
We develop a general and consistent finite approximation scheme with tractable solutions and error bounds for computing the stationary distribution and associated steady-state performance measures of an MC supported on $\mathbb{R}$ (Theorem \ref{Theorem: error bound naive},\ref{theorem:performance measure error},\ref{Theorem: convergence naive}). 
The MC's support can be discrete, continuous or even mixed.
Our approximate MC is constructed based on variation properties of the original MC's transition kernel, and the error is obtained from a linear program with a relatively low computation cost. 
%All the above features are useful when we apply our approach to decision problems.

%; see examples in \cite{li2021modeling}. 
%Therefore, our finite approximation approach provides a significant advance in stochastic modeling and analysis.

Compared with simulation approaches using an empirical distribution as a \textit{proxy stationary distribution} (e.g., Markov chain Monte Carlo (MCMC), discrete event simulation, etc.), our approach has better tractability because our approximate stationary distribution obtained from the constructed finite-state MC is explicit. 
The nature of error in these approaches is different. 
For simulation, the error of empirical distribution is usually quantified as probabilistic concentration bounds. 
In contrast, the finite approximation error is deterministic and guaranteed in the worst case. 

%\noindent \textbf{Related applications.}\quad
%Although this paper focuses on discrete time MC, the results are applicable to analysis of continuous time Markov processes (e.g., Brownian motions, Feller processes) via their embedded discrete time MC; see examples of using embedded MC in \cite{bremaud2013markov}. In addition to stochastic models' performance evaluation and optimization problems, our MC finite approximation approach is also useful for approximate Markov chain Monte Carlo (MCMC) in the analysis of big data and Bayesian inference: making approximations to the Markov transition kernel is a common technique to control computation time, and our results indicate how transition kernel differences of MC are reflected in differences of their distributions. Recent works on approximate MCMC techniques include \cite{pillai2014ergodicity,johndrow2015optimal,rudolf2018perturbation}. 

\noindent \textbf{An Illustration Example of G/G/1+G Queue}. \quad Our method can be used for any stochastic system that can be represented by an MC satisfying the settings in this paper. 
As an example, we consider a G/G/1+G queue, where a customer abandons the system once her deadline elapses before service is offered. 
The distributions of inter-arrival times, service times and patience times are all general, and a first-come-first-serve scheme is assumed. 
To the best of our knowledge, there is no known analytical or numerical solution for computing the \textit{virtual waiting time}'s stationary distribution in such a queue.
We compare our approach against benchmark results including an MCMC implementation, and a fluid approximation \cite{whitt2006fluid}. We find that our approach is significantly more accurate and efficient in obtaining the stationary distribution and evaluating a few illustrative example steady-state performance measures including abandonment probability, ``no wait" probability, queue length, etc.%; Probability of no wait, Probability of Abandonment, Expected Queue Length.
%Additional applications and computation experiments with other stochastic models can be found in our companion paper \cite{li2021modeling}. 

\noindent \textbf{Efficient and Near-Optimal Approximation.}\quad 
Under a Lipschitz continuity assumption on the original MC's transition kernel, our approach's supremum-normed error in computing stationary distributions can be as small as $O({N}^{-1})$, where $N$ is the number of jumps (or states) used to construct the approximate MC.
Moreover, our approach is near-optimal among all approximation methods using discrete distributions including the empirical distributions obtained by a simulation approach (Theorem \ref{theorem:near-optimality}).  
The empirical distributions can also be regarded as a discrete proxy stationary distribution for the original MC but converge as $O({N}^{-\frac{1}{2}})$ according to the Dvoretzky–Kiefer–Wolfowitz (DKW) concentration inequality (Theorem 11.5 in \cite{kosorok2006introduction}), where $N$ is the number of jumps (or sample size) used in the empirical distribution. 
%The efficiency of finite approximation compared against computer simulation is confirmed via numerical experiments in Section \ref{sec:Numerical}. 

\noindent \textbf{Weak Assumptions.}\quad Our results are developed under the supremum norm, which has several advantages over the conventional total variation norm. 
The supremum norm requires weaker assumptions on guaranteeing solution consistency.
Moreover, it allows us to capture many consistency situations that are invalid under the total variation norm: when we apply finite approximation and use a finite-state MC to approximate another MC with a continuous stationary distribution, the total variation-normed error always equals $1$.
Thus, the total variation norm is not an adequate measure. 
In contrast, the supremum-normed error can approach $0$ in a proper approximation. % (see examples in Section \ref{sec:Numerical}). 

Similar to \cite{herve2014approximating}, we do not require approximate MCs' transition kernels converging under the operator norm, which is often seen in perturbation analysis but can be a strong assumption in practice. 
Our consistency results require compactness and stationary distribution uniqueness of the original MC's transition kernel, while our error bound results require a weaker condition than the approximate finite-state MC's irreduciblity. 
Ergodicity %is not necessary but it may
can serve as a sufficient condition for computing error bounds (see Section \ref{sec:error bound naive}). 

\noindent \textbf{A Novel Integral Equation Analysis Framework.}\quad
Our results are based on a key observation that by adding a natural scaling constraint, many balance equations supported on $\mathbb{R}$ can be transformed into a canonical-form distributional integral equation. 
%The approximate solution's error bound is developed by considering the equation in a specially designed Banach space of distributions, and the consistency theory is developed by considering the transition kernel as a compact operator whose finite-rank approximations are consistent. 
This class of distributional integral equations and their numerical solutions were analyzed in our companion paper \cite{li2021numerical}, which provide a theoretical foundation for the current work. %This paper further develops the results of \cite{li2021numerical} in the context of Markov Chains.

\subsection{Related Literature.}   \label{sec:LiteratureReview} \,  \\
\noindent \textbf{MC Finite Approximation.}\quad 
A significant literature has investigated the approximation of a countable-state MC by another (Table \ref{table:lit rev}).
%Methods include truncation-augmentation \citep{seneta1967finite, heyman1991approximating, zhao1996censored, tweedie1998truncation, hart2012convergence, liu2015perturbation, liu2018error2, liu2018error,masuyama2016error, masuyama2016limit, masuyama2017continuous, masuyama2017error}, block-based method \citep{bright1995calculating, li2000stochastic, baumann2010numerical, baumann2013computing, phung2010simple, dayar2011infinite,dayar2012kronecker}, aggregation-disaggregation \citep{haviv1987aggregation},  nearly complete decomposition \citep{haviv1986approximation}, etc. 
%A more general method for countable-state MC  is given in \cite{wolf1980approximation}. 
However, these methods rely on the transition matrix structure, and are not extended for continuous-state MCs.
\cite{pincus1992approximating} used finite-state MCs to approximate a continuous-state MC under a strong assumption of kernel continuity and provided a weak convergence result of stationary distributions. 
In practice, however, jumps and kernel discontinuities are common, and a weak convergence without tractable error bounds has limited applications. 
More recently, \cite{herve2014approximating, herve2019state, herve2020v} used finite-rank sub-Markov kernels (which may not be a transition kernel of an MC) to approximate $V$-geometric ergodic MC kernels that could be continuous-state. 
However, $V$-geometric ergodicity can be a strong assumption \citep{roberts2008variance}. 
In addition, the authors measured approximation errors by the total variation norm. 
As discussed in subsection \ref{sec:contributions}, results using the total variation norm are weaker than that based on the supremum norm. 

\begin{table}[!htp]
\scriptsize
\caption{Summary of Mrkov chain finite approximation approaches in literature.}
\label{table:lit rev}
\centering
\begin{tabular}{|l|l|l|l|l|l|}
\hline
\multicolumn{2}{|l|}{Classification}                                                                                                                & Model                                                     & Key Assumptions                                                                                                                              & Main Results                                                                                                              & References                                            \\ \hline
\multicolumn{2}{|l|}{\begin{tabular}[c]{@{}l@{}}Computational  algorithms\end{tabular}}                                                                & \begin{tabular}[c]{@{}l@{}}Countable \\ state\end{tabular} & \begin{tabular}[c]{@{}l@{}}Ergodic level dependent quasi-birth- \\and-death process. \end{tabular}                                             & \begin{tabular}[c]{@{}l@{}}Efficient algorithms, accuracy \\  validated by  experiments.\end{tabular}                       & \begin{tabular}[c]{@{}l@{}}\cite{baumann2010numerical,phung2010simple,dayar2011infinite} \, \end{tabular}     \\ \hline
\multicolumn{2}{|l|}{\begin{tabular}[c]{@{}l@{}}Perturbation  analysis\end{tabular}}                                                               & \begin{tabular}[c]{@{}l@{}}Countable \\ state\end{tabular} & \begin{tabular}[c]{@{}l@{}}Irreducibility,  positive recurrence, \\  geometric ergodicity, uniform \\  ergodicity, etc.\end{tabular}             & Explicit error bounds.                                                                                                   & \begin{tabular}[c]{@{}l@{}}\cite{tweedie1980perturbations,roberts1998convergence,altman2004perturbation,hunter2005stationary}\\ \cite{liu2012perturbation,mouhoubi2010new}\end{tabular} \\ \hline
\multirow{2}{*}{\begin{tabular}[c]{@{}l@{}}Truncation \\ methods\end{tabular}} & \begin{tabular}[c]{@{}l@{}}Convergence \\ analysis\end{tabular}    & \begin{tabular}[c]{@{}l@{}}Countable \\ state\end{tabular} & \begin{tabular}[c]{@{}l@{}}Irreducibility, positive recurrence,\\ stochastic  monotonicity, etc.\end{tabular}                                 & \begin{tabular}[c]{@{}l@{}}Convergence in $L^1$ norm, $V$\\ norm, total  variation norm, etc.\end{tabular}                 & \begin{tabular}[c]{@{}l@{}}\cite{seneta1967finite,seneta1980computing, gibson1987augmented,gibson1987monotone}\\ \cite{heyman1991approximating,hart2012convergence,liu2015perturbation} \end{tabular} \\ \cline{2-6} 
                                                                               & \begin{tabular}[c]{@{}l@{}}Error \\ analysis\end{tabular}          & \begin{tabular}[c]{@{}l@{}}Countable \\ state\end{tabular} & \begin{tabular}[c]{@{}l@{}}Irreduciblility, positive recurrence,\\ upper Hessenber structure,\\ modulated drift condition, etc.\end{tabular} & \begin{tabular}[c]{@{}l@{}}Errors in $L^1$ norm,  supremum \\ norm, $V$ norm or total \\ variation  norm.\end{tabular} & \begin{tabular}[c]{@{}l@{}}\cite{zhao1996censored,tweedie1998truncation,liu2010augmented,liu2018error2}\\\cite{liu2018error} \, \end{tabular} \\ \hline
\multicolumn{2}{|l|}{\begin{tabular}[c]{@{}l@{}}Block-based methods\end{tabular}}                                                                & \begin{tabular}[c]{@{}l@{}}Countable \\ state\end{tabular} & \begin{tabular}[c]{@{}l@{}}Block-structure, and block-\\monotonicity.\end{tabular}                                                          & \begin{tabular}[c]{@{}l@{}}Pointwise convergence,  errors \\ in  supremum norm or  total \\ variation  norm.\end{tabular}       & \begin{tabular}[c]{@{}l@{}} \cite{li2000stochastic,masuyama2015error,masuyama2016error,masuyama2017continuous} \, \end{tabular} \\ \hline
\multirow{3}{*}{\begin{tabular}[c]{@{}l@{}}General\\ methods\end{tabular}}     & \begin{tabular}[c]{@{}l@{}}Countable \\ state problem\end{tabular} & \begin{tabular}[c]{@{}l@{}}Countable \\ state\end{tabular} & \begin{tabular}[c]{@{}l@{}}Irreduciblility, and positive\\  recurrence.\end{tabular}                                                          & Pointwise convergence.                                                                                                   & \begin{tabular}[c]{@{}l@{}}   \cite{wolf1980approximation}   \, \end{tabular}                                             \\ \cline{2-6}     
     & \begin{tabular}[c]{@{}l@{}}Continuous \\ state problem\end{tabular}   & \begin{tabular}[c]{@{}l@{}}Continuous \\ state\end{tabular}  & \begin{tabular}[c]{@{}l@{}}Transition kernel's weak continuity \\  and uniform  weak convergence.\end{tabular}                 & \begin{tabular}[c]{@{}l@{}} Convergence in  distribution.\end{tabular}                               & \begin{tabular}[c]{@{}l@{}} \cite{pincus1992approximating} \end{tabular}     \\ \cline{2-6}
                                                                               & \begin{tabular}[c]{@{}l@{}}General \\ state problem\end{tabular}   & \begin{tabular}[c]{@{}l@{}}General \\ state\end{tabular}  & \begin{tabular}[c]{@{}l@{}}$V$-geometric ergodicity,  weakened\\ convergence,  and uniform weak drift.\end{tabular}                 & \begin{tabular}[c]{@{}l@{}}Convergence and errors in total \\ variation norm.\end{tabular}                               & \begin{tabular}[c]{@{}l@{}} \cite{herve2014approximating,herve2019state,herve2020v} \end{tabular}     \\ \hline
\end{tabular}
\end{table}

\noindent \textbf{Comparison with MCMC.}\quad
Both MCMC and finite approximation evaluate the steady-state performance via generating a proxy stationary distribution. 
For the finite approximation, the proxy stationary distribution is the approximate finite-state MC's stationary distribution, while for MCMC, the proxy stationary distribution is the empirical distribution of the sample.
For MCMC, convergence rates and concentration bounds have been extensively studied, e.g., \cite{gelman1992inference, liu1993control,  yu1995estimating, brooks1998convergence, brooks1998general, kontoyiannis2012geometric, bremaud2013markov}. For a more systematic review, we refer the readers to \cite{han2001markov, sinharay2003assessing, brooks2011handbook}. %A general method used to analyze a geometrically ergodic chain is via the spectrum analysis of its transition kernel. Results on convergence rates and error bounds are obtained in the Hilbert space $L^2 = L^2(p)$, where $p$ denotes the MC  's stationary  distribution. For example, a well-known result from  \cite{roberts1997geometric} states that for a reversible, $\Psi$-irreducible and aperiodic discrete time MC   with initial distribution $\mu$, we have 

\subsection{Organization.}
We close this section with an outline of this paper. 
Section \ref{sec:prelim} introduces key notations and settings. 
Section \ref{sec:error bound naive} analyzes our finite approximation approach that uses an approximate finite-state MC to compute the orignal MC's stationary distribution. 
We provide error bounds for the approximate stationary distribution as well as the corresponding steady-state performance measures. 
Section \ref{sec:convergence naive} shows that under certain regularity conditions, the stationary distribution and steady-state performance measures computed by our approach are consistent under the supremum norm. 
Section \ref{sec:verify RC} provides verifiable sufficient conditions ensuring the regularity conditions defined in Section \ref{sec:convergence naive}. 
Section \ref{sec:Numerical} illustrates the use of our approach with an example of G/G/1+G queue, and compares against benchmark results. 
Section \ref{sec:ConcludingRemarks} presents concluding remarks. 
All cited theorems from the literature are listed in Appendix \ref{append:cited theorems}.   

\section{Preliminaries.}\label{sec:prelim} 
We now present definitions and assumptions for analyzing our finite approximation. 
We assume support $\Omega=[0,1]$ for the general results of this paper. 
Our analysis can be adapted for other support $\Omega'\subseteq \mathbb{R}$,  e.g., transforming an unbounded support 
$\mathbb{R}_+$ into $(0,1]$ via function $e^{-x}$ and then expanding support $(0,1]$ to $[0,1]$ via adding a (trivial) state $0$. 
Our results (Theorem \ref{Theorem: error bound naive},\ref{theorem:performance measure error},\ref{Theorem: convergence naive}), conditions (Condition \ref{A_UBV}-\ref{A_T2}) and corresponding proofs can also be similarly developed for $\Omega'$.

\subsection{Notations and a Banach Space.}\label{sec:banach space}
We first construct a Banach space of right-continuous distribution functions. 
Let $\mathbf{D}$ be the collection of probability distribution functions of random variables defined on $(\Omega,\mathcal{B})$ with $\mathcal{B}$ being the Borel algebra for $\Omega$. 
Define $\mathbf{X}:=$\rm span$(\mathbf{D})= \left\{\sum_{k=1}^{n} a_k f_k\,|\, n\in \mathbb{N}_{++}, a_k \in \mathbb{R}, f_k \in \mathbf{D} \right\}$, which is also the linear space of distribution functions of all finite signed measures on $(\Omega,\mathcal{B})$. 
Let $||f||_\infty$ be the supremum norm of $f:\mathbb{R}\mapsto\mathbb{R}$, and $\bar{\mathbf{X}}$ be the closure of $\mathbf{X}$ with norm $||\cdot||_\infty$. 
Then, $(\bar{\mathbf{X}}, ||\cdot||_\infty)$ is a Banach space (Theorem 2.5 in \cite{li2021numerical}). 

We next define operations on space $\bar{\mathbf{X}}$.
We write ``$f_k$ converging to $f$ on $\bar{\mathbf{X}}$" as ``$f_k\rightrightarrows f$" to emphasize uniform convergence. Let $V(f;\Gamma)$ be the total variation of $f:\mathbb{R}\mapsto \mathbb{R}$ on $\Gamma\subseteq \mathbb{R}$:
\begin{align}
    V(f;\Gamma):= \sup_{x_1<x_2<...<x_n; x_1,x_2,...,x_n\in\Gamma;n\in\mathbb{N} } \sum_{i=2}^n|f(x_i)-f(x_{i-1})|, \nonumber
\end{align}
where the supremum operation is over all partitions of $\Gamma$. 
For convenience, let $V(f):= V(f;\Omega)$. 
For function $f(x,u):\mathbb{R}^n\times \mathbb{R}\mapsto\mathbb{R}$, we use $V_u(f(x,u))$ to denote the total variation of $f(x,u)$ as a single-variable function of $u$ with any fixed $x\in \mathbb{R}^n$.  
%The total variation $V(f)$ measures the one-dimensional length of the function curve of $f$.  
Let $\mathcal{I}$ be the identity operator and $||\cdot||_O$ be the operator norm for an operator $\mathcal{L}$ on $\bar{\mathbf{X}}$: $||\mathcal{L}||_O = \sup_{f\in \bar{\mathbf{X}}, ||f||_\infty = 1} ||\mathcal{L}f||_\infty.$

We now define \textit{finite (transition) kernels} $T(\Omega, \mathcal{B})$ in the distribution sense.
\begin{definition}[Finite Kernel] \label{Def:Kernel}
Let $T(\Omega, \mathcal{B})$ be the collection of $\tau(x,u):\Omega^2\mapsto \mathbb{R}_+$ such that $(i)$ for all $x\in\Omega$, $\tau(x,\cdot)$  is $\mathcal{B}$-measurable, and $(ii)$ for all $u\in\Omega$, $\tau(\cdot,u)=0$ or there exists $\alpha>0$ such that $\alpha  \tau(\cdot,u)\in\mathbf{D}$. Particularly, $\tau\in T(\Omega, \mathcal{B})$ is a Markov kernel if $\tau(\cdot,u)\in\mathbf{D}$ for all $u\in\Omega$.
\end{definition}
Both finite and Markov kernels will be used to describe MCs' transitions and balance equations.
It is possible to formulate a kernel $\tau\in {T}(\Omega, \mathcal{B})$  into a continuous linear operator on space $\bar{\mathbf{X}}$.
\begin{definition}[Transition Operator] 
\label{Def:K}
For kernel $\tau \in {T}(\Omega,\mathcal{B})$, its corresponding transition operator $\mathcal{L}$ on $\bar{{\mathbf{X}}}$ is defined as follows: if $f\in \mathbf{X}$, then  
\begin{align}
    \mathcal{L} f(x) := \int_\Omega \tau(x,u)df(u),\quad \forall x\in \Omega. \label{Def:K_1}
\end{align}
Otherwise, we select an arbitrary sequence $\{f_k\}_{k\in \mathbb{N}}\subseteq \mathbf{X}$ that converges to $f$ on $\mathbf{\bar{X}}$ and define
\begin{align}
    \mathcal{L} f(x) := \lim_{k\rightarrow \infty}\mathcal{L} f_k(x),\quad \forall x\in \Omega. \label{Def:K_2}
\end{align}
\end{definition}
Then we can use $f'=\mathcal{L}f$ to describe that an initial probability distribution $f$ returns distribution $f'$ after one transition in an MC with kernel $\tau$.
According to Theorem 3.1 of \cite{li2021numerical}, transition operator $\mathcal{L}$ in Definition \ref{Def:K} is a well-defined continuous linear operator iff the following holds. 
\begin{condition}
[Uniformly Bounded Variation] \label{A_UBV} For $\tau\in{T}(\Omega, \mathcal{B})$, $ \sup_{x\in\Omega}{V_u(\tau(x,u))}<\infty$.
\end{condition}
According to Condition \ref{A_UBV}, for all $x\in \Omega$, $\tau(x,u)$ as a single-variable function of $u$ has a uniformly bounded total variation.
In Definition \ref{Def:K}, $\mathcal{L}f$ is not defined uniformly, i.e., \eqref{Def:K_1} vs. \eqref{Def:K_2}. This is because  function $f\in \mathbf{\bar{X}}\backslash\mathbf{X}$ as a distribution of measures is ill defined and the conventional definition \eqref{Def:K_1} does not apply. 
Instead, we define \eqref{Def:K_2} via convergence, which is natural because any continuous operator must have \eqref{Def:K_2} as a property. 
With slight abuse of notation, we also write $\mathcal{L} f(x)$ as $\int_\Omega \tau(x,u)df(u)$ for all $f\in\mathbf{\bar{X}}\backslash \mathbf{X}$. %The set $\mathbf{\bar{X}}\backslash\mathbf{X}$ does not include our original or approximate stationary distributions and will  be mainly used for computing error bounds and proving convergence of approximate stationary distributions.

\subsection{A Finite Approximation Approach.} \label{sec:regularity conditions}
Our approximation approach is given below.
\begin{approximation}\label{A:Stationary Dist Eq}
Consider a discrete-time MC supported on $\Omega$, denoted as the original MC. Its stationary distribution $p\in\mathbf{D}$ satisfies the following balance equations (by letting $f=p$)
\begin{align}
    f(x) =\int_{\Omega} \bar{\kappa}(x,u) df(u), \quad \forall x \in \Omega, \label{eqn:MC_BE_original} \tag{BE-1}
\end{align}
where the Markov kernel $\bar{\kappa}\in T(\Omega, \mathcal{B})$. 
Also consider a sequence of discrete-time finite-state MCs supported on $\Omega$, denoted as the approximate MCs and indexed by $r\in\mathbb{N}$. 
For the $r$-th approximate MC, the states are  $0=c_1^{(r)}<c_2^{(r)}<...<c_{J^{(r)}}^{(r)}=1$, the transition matrix is  $Q^{(r)}=\{q^{(r)}_{ij}\}_{i,j=1,2,...,J^{(r)}}$, and the stationary distribution is $p^{(r)}$.
We compute the original MC's steady-state performance measure $\int_{\Omega} \phi(x)dp(x)$ by using $p\approx p^{(r)}.$ 
\end{approximation}
Section \ref{sec:convergence naive} will show how to properly construct the approximate MCs. % so that the approximate stationary distributions $\{p^{(r)}\}_{r\in\mathbb{N}}$ converge to the original stationary distribution $p$ under the supremum norm. 
We describe the $r$-th approximate MC's balance equations in a distributional form, similar to \eqref{eqn:MC_BE_original} for the original MC:
\begin{align}
    f(x) = &   \int_\Omega \bar{\kappa}^{(r)}(x,u) df(u),\quad \forall x \in \Omega.\label{eqn:MC_BE_approx} \tag{ABE-1}
\end{align}
Here transition kernel $\bar{\kappa}^{(r)}$ is specified by the approximate finite-state MC's transition matrix $Q^{(r)}$:
\begin{align}
    &\bar{\kappa}^{(r)}(x,u) :=  \sum_{i=1}^{J^{(r)}}\bar{\omega}_i^{(r)}(x)\mathbf{1}\{u\in(c^{(r)}_{i-1},c^{(r)}_i]\},\quad  \forall  x,u\in\Omega,r\in\mathbb{N},  \label{eqn:def kappa}\\
    &\bar{\omega}_i^{(r)}(x) := \sum_{j=1}^{J^{(r)}}\mathbf{1}\{c_j^{(r)}\leqslant x\}q^{(r)}_{ij},\quad x\in\Omega,i=1,2,...,J^{(r)},r\in\mathbb{N}, \label{eqn:def bar omega} 
\end{align}
and $c_0^{(r)}:= -\infty$. 
%It is not hard to verify that $\bar{\omega}_{i}^{(r)}\in\bar{\mathbf{X}}\,(i=1,2,...,J^{(r)})$ and $\bar{\kappa}^{(r)}\in T(\Omega, \mathcal{B})$ satisfies Condition \ref{A_UBV} $(r\in\mathbb{N})$.
Then $p^{(r)}$ is a solution to \eqref{eqn:MC_BE_approx} (by letting $f =p^{(r)}$, $r\in\mathbb{N}$).

\section{Errors of Finite Approximation. }\label{sec:error bound naive}
%We first consider the error in using the stationary distribution of a finite-state MC to approximate that of the original MC. Our main result %in the following theorem 
Our main result of this section is the following bound.
\begin{theorem}[Errors in Stationary Distributions]\label{Theorem: error bound naive}
Consider an instance of approximate MC in Approximation \ref{A:Stationary Dist Eq}.  
If the original MC's kernel $\bar{\kappa}$ satisfies Condition \ref{A_UBV} and the approximate MC is irreducible, then the error $||p-p^{(r)}||_\infty \leqslant e_1\cdot e_2$, where $e_1= \sup_{x,u\in\Omega} |\bar{\kappa}(x,u) -\bar{\kappa}^{(r)}(x,u)|$, $e_2^{-1}=\min_{k=0,1,...,J^{(r)}}\{y^*_{k}\}$, and $y^*_{k}$ is obtained from a linear program with $a_0\equiv 0$ $(k=0,1,...,J^{(r)})$:
\begin{align}
    &\min_{y_{k}\in\mathbb{R},\{a_j\}_{j=1}^{J^{(r)}} \in[-1,1]^{J^{(r)}} } \quad\, y_{k}  \nonumber\\
    &s.t.  \quad \,  
    -y_{k}  \leqslant a_j+ \eta \delta_{jk}(1-a_j) -\sum_{i=1}^{J^{(r)}} (1+\sum_{s=1}^{j}q^{(r)}_{is})\cdot(a_i-a_{i-1}) \leqslant y_{k},\,\quad j=0,1,...,J^{(r)}, \eta\in \{0,1\}.  \nonumber
\end{align}   
%{\bf Possibly give the Linear Programs here, and delete Lemma 3. Also give a formal proof by invoking the result in Li and Mehrotra, Lemma 1 and Lemma 2}
\end{theorem}

Theorem \ref{Theorem: error bound naive} provides the following insight.
%follows classical algebraic results for linear equation system sensitivity analysis: 
The balance equations of both the original and approximate MCs are linear systems. 
The later can be regarded as a computable approximation of the former with a perturbation.
Our results indicate that the approximation error, i.e., the solution error due to perturbation, is decided by a factor $e_1$ measuring the MC kernel difference (original vs. approximate), and factor $e_2$ measuring the approximate stationary distribution's sensitivity.

Our proof is completed in three steps. 
First, we show that the original and approximate MCs' balance equations can both be formulated into operator-form equations  (Subsection \ref{sec:IE reform}). 
Next, we obtain an operator-form approximation error bound (Subsection \ref{sec:operator form error}).
Lastly, we show the operator-form bound can be computed by $e_1\cdot  e_2$ (Subsection \ref{sec:compuate error naive}). 

The error bounds for steady-state performance measures are provided in Subsection \ref{sec:performance error}.

\subsection{Transformation of Balance Equations.}\label{sec:IE reform}
We show that by adding a natural scaling constraint, balance equations of both the original and approximate MCs can be transformed into operator equations. 
Our transformation brings a key advantage: balance equations \eqref{eqn:MC_BE_original} or \eqref{eqn:MC_BE_approx} have at least one degree of freedom in its solution space, while by incorporating a natural scaling constraint, the transformed equations have a unique solution under appropriate assumptions. 
Our error bound and consistency analysis will both be based on this solution uniqueness.
\begin{lemma}[Balance Equation Reformulation I]\label{Lemma:IE reform I}
%Consider an original MC in Definition \ref{A:Stationary Dist Eq} with stationary distribution $p$ and kernel $\bar{\kappa}$ satisfying Condition \ref{A_UBV}, and one instance of approximate finite-state MC in Definition \ref{A:Stationary Dist Eq Approx} with index $r\in\mathbb{N}$ and stationary distribution $p^{(r)}$. 
Given the settings and assumptions in Theorem \ref{Theorem: error bound naive}, $f=p$ and $f=p^{(r)}$ are respectively solutions to %the following integral equations
\begin{align}
     & f(x)    =  \xi(x)   + \int_{\Omega} \kappa(x,u) d f(u), \quad \forall x\in \Omega; \label{eqn:MC_IE_original} \tag{BE-2} \\
     & f(x)  =   \xi(x)  +  \int_{\Omega} \kappa^{(r)}(x,u) d f(u),\quad \forall x\in \Omega, \label{eqn:MC_IE_approx} \tag{ABE-2}
\end{align}
where the inhomogeneous term $\xi\in\mathbf{\bar{X}}$ and the kernels $\kappa,\kappa^{(r)}\in{T}(\Omega, \mathcal{B})$  are specified as: $\forall x,u\in\Omega,$
\begin{align*}
    & \xi(x):=-\mathbf{1}\{x\geqslant 0\}, \qquad\,\,\,\,\,  \kappa(x,u) := \mathbf{1}\{x\geqslant 0\}+\bar{\kappa}(x,u),  \\
    & \kappa^{(r)}(x,u) :=  \sum_{i=1}^{J^{(r)}}[\mathbf{1}\{x\geqslant 0\}+\bar{\omega}_i^{(r)}(x)]\cdot \mathbf{1}\{u\in(c^{(r)}_{i-1},c^{(r)}_i]\}. 
    %\label{eqn:def xi}  \label{eqn:def kappa} \label{eqn:def kappa r}
\end{align*}
Equivalently \eqref{eqn:MC_IE_original} and \eqref{eqn:MC_IE_approx} can be stated as %an operator equation on space $\bar{\mathbf{X}}$:
\begin{align}
    & f = \xi+\mathcal{K}f, \label{eqn:MC_OE_original}\tag{BE-3}\\
    & f =\xi + {\mathcal{K}^{(r)}} f, \label{eqn:MC_OE_approx} \tag{ABE-3}
\end{align}
where operators $\mathcal{K}$ and $\mathcal{K}^{(r)}$ are respectively defined by Definition \ref{Def:K} with kernels $\kappa$ and $\kappa^{(r)}$.
\end{lemma}
\begin{proof}{Proof:} %We first show the transformation into \eqref{eqn:MC_IE_original} and \eqref{eqn:MC_IE_approx}. 
%Consider the original MC. 
Recall that $p$ satisfies \eqref{eqn:MC_BE_original} and an extra scaling constraint $p(1)=1$ because $p$ is a probability distribution supported on $\Omega=[0,1]$. 
Thus, $f=p$ satisfies each of the following systems. 
\begin{align}
    & f(1)=1;\,f(x)= \int_{\Omega} \bar{\kappa}(x,u) df(u), \, \forall x \in \Omega. \nonumber  \\
    \Rightarrow & f(x)= \mathbf{1}\{x\geqslant 0\} \cdot [f(1)-1] +  \int_{\Omega}  \bar{\kappa}(x,u)df(u),\,  \forall x\in\Omega. \qquad \text{(due to $\mathbf{1}\{x\geqslant 0\} \cdot [f(1)-1]=0$)} \nonumber \\
    \Rightarrow & f(x)= -\mathbf{1}\{x\geqslant 0\}   +  \int_{\Omega} \big[\mathbf{1}\{x\geqslant 0\}+\bar{\kappa}(x,u) \big] df(u),\,  \forall x\in\Omega. \qquad \text{(due to $f(1)=\int_\Omega df(u) $)} \nonumber\\
    \Rightarrow & f(x)= \xi(x)  +  \int_{\Omega} \kappa(x,u) df(u),\,  \forall x\in\Omega.\qquad \text{(due to definitions of $\xi$ and $\kappa$)}  \label{IE reform 1: right}
\end{align} 
Thus, $f=p$ is a solution to \eqref{eqn:MC_IE_original}. 
We next consider %the approximate finite-state MC. Recall that stationary distribution 
$p^{(r)}$, which similarly satisfies \eqref{eqn:MC_BE_approx} and an extra scaling constraint $p^{(r)}(1)=1$ because $p^{(r)}$ is a probability distribution supported on $\Omega=[0,1]$. Thus, $f=p^{(r)}$ satisfies each of the following systems. 
\begin{align}
    & f(1)=1;\,f(x)= \int_{\Omega} \bar{\kappa}^{(r)}(x,u) df(u), \, \forall x \in \Omega. \nonumber \\
    %\Rightarrow & f(c^{(r)}_{J^{(r)}})=1;\,f(x)= \int_{\Omega} \bar{\kappa}^{(r)}(x,u) df(u), \, \forall x \in \Omega. \nonumber \\
    \Rightarrow & f(x)= \mathbf{1}\{x \geqslant 0\} \cdot [f(1)-1] +  \int_{\Omega}  \bar{\kappa}^{(r)}(x,u)df(u),\,  \forall x\in\Omega.  \nonumber \qquad \text{(due to $\mathbf{1}\{x\geqslant 0\} \cdot [f(1)-1]=0$)} \\
    \Rightarrow & f(x)= -\mathbf{1}\{x \geqslant 0\}   +  \int_{\Omega} \big[\mathbf{1}\{x \geqslant 0\}+\bar{\kappa}^{(r)}(x,u) \big] df(u),\,  \forall x\in\Omega. \nonumber \qquad \text{(due to $f(1)=\int_\Omega df(u) $)} \\
    \Rightarrow & f(x)= \xi(x)  +  \int_{\Omega} \sum_{i=1}^{J^{(r)}}[\mathbf{1}\{x \geqslant 0\}+\bar{\omega}_i^{(r)}(x)]\cdot \mathbf{1}\{u\in(c^{(r)}_{i-1},c^{(r)}_i]\} df(u),\,  \forall x\in\Omega. \nonumber \\
    & \text{(due to the definition of $\bar{\kappa}$ and that $\Omega=\cup_{i=1}^{J^{(r)}}\{(c^{(r)}_{i-1},c^{(r)}_i]\cap\Omega\}$)} \nonumber \\
    \Rightarrow & f(x)= \xi(x)  +  \int_{\Omega} \kappa^{(r)}(x,u) df(u),\,  \forall x\in\Omega. \nonumber \qquad \text{(due to definitions of $\xi$ and $\kappa^{(r)}$)} 
\end{align} 
Thus,  $f=p^{(r)}$ is a solution to \eqref{eqn:MC_IE_approx}. %Once we transform the original and approximate MC's balance equations into \eqref{eqn:MC_IE_original} and \eqref{eqn:MC_IE_approx}, 
Lastly, we directly have the operator forms \eqref{eqn:MC_OE_original} for \eqref{eqn:MC_IE_original} and \eqref{eqn:MC_OE_approx} for \eqref{eqn:MC_IE_approx} via Definition \ref{Def:K}. 
% Note that because $\bar{\kappa}$ satisfies Condition \ref{A_UBV},  $\kappa= \bar{\kappa} + \mathbf{1}\{x \geqslant 0\}$ also satisfies Condition \ref{A_UBV}. It is not hard to verify that $\kappa^{(r)}$ also satisfies Condition \ref{A_UBV} due to its finite support. Therefore,  operators $\mathcal{K}$ and $\mathcal{K}^{(r)}$ are both well-defined (Theorem 3.1 in \cite{li2021numerical}). Then  \eqref{eqn:MC_OE_original} and \eqref{eqn:MC_OE_approx} are also well-defined equations.
% This completes the proof of Lemma \ref{Lemma:IE reform I}. 
\Halmos
\end{proof}

\subsection{Operator-Form Error Bound.} \label{sec:operator form error}
Lemma \ref{Lemma:IE reform I} implies the following error bound. 

\begin{lemma}[Operator-Form Errors in Stationary Distributions]\label{Lemma:operator form error naive}
%Consider an original MC in Definition \ref{A:Stationary Dist Eq} with stationary distribution $p$ and kernel $\bar{\kappa}$ satisfying Condition \ref{A_UBV}, and one instance ofAccording to Lemma \ref{Lemma:IE reform I}, error factor $||(\mathcal{I}-\mathcal{K}^{(r)})^{-1}||_O$ in \eqref{eqn:error bound naive} measures the approximate stationary distribution's sensitivity (in the inhomogeneous term $\xi$) while factor $||\mathcal{K}p- \mathcal{K}^{(r)}p ||_\infty$ measures the kernel difference of the original vs.  approximate MC. approximate finite-state MC in Definition \ref{A:Stationary Dist Eq Approx} with index $r\in\mathbb{N}$ and stationary distribution $p^{(r)}$.
Given the settings and assumptions in Theorem \ref{Theorem: error bound naive}, we have the following bound
\begin{align}
    \left|\left| p-  p^{(r)} \right|\right|_\infty \leqslant ||(\mathcal{I}-\mathcal{K}^{(r)})^{-1}||_O \cdot ||\mathcal{K}p- \mathcal{K}^{(r)}p ||_\infty.
    \label{eqn:error bound naive}
\end{align}
%as long as the inverse $(\mathcal{I}-\mathcal{K}^{(r)})^{-1}$ exists.
%where operators $\mathcal{K}$ and $\mathcal{K}^{(r)}$ are respectively defined by kernels $\kappa$ and $\kappa^{(r)}$ with Lemma \ref{Lemma:IE reform I}. 
\end{lemma} 

\begin{proof}{Proof:}
Suppose $(\mathcal{I}-\mathcal{K}^{(r)})^{-1}$ exists. 
Because $ p = \xi+\mathcal{K}p$ and $p^{(r)} =\xi + {\mathcal{K}^{(r)}} p^{(r)}$ by Lemma \ref{Lemma:IE reform I},
\begin{align*}
    \left|\left|p-p^{(r)}\right|\right|_\infty =&\left|\left|(\mathcal{I}-{\mathcal{K}}^{(r)})^{-1}\left( p-{\mathcal{K}}^{(r)}p\right)-(\mathcal{I}-{\mathcal{K}}^{(r)})^{-1}\xi \right|\right|_\infty=\left|\left|(\mathcal{I}-{\mathcal{K}}^{(r)})^{-1}\left( p-\xi-{\mathcal{K}^{(r)}}p\right)\right|\right|_\infty \\
    \leqslant & \left|\left|(\mathcal{I}-{\mathcal{K}^{(r)}})^{-1}\left( {\mathcal{K}}p-{\mathcal{K}^{(r)}}p\right)\right|\right|_\infty 
    \leqslant  ||(\mathcal{I}-{\mathcal{K}}^{(r)})^{-1}||_O\cdot\left|\left| {\mathcal{K}}p-{\mathcal{K}^{(r)}}p\right|\right|_\infty.
\end{align*}
which is the bound in \eqref{eqn:error bound naive}.
Next we show that $(\mathcal{I}-\mathcal{K}^{(r)})^{-1}$ indeed exists. Define weight functions 
\begin{align}
    {\omega}_i^{(r)}(x) := & \mathbf{1}\{x\geqslant 0\}+\bar{\omega}_{i}^{(r)}(x)= \mathbf{1}\{x\geqslant 0\}+\sum_{j=1}^{J^{(r)}}\mathbf{1}\{c_j^{(r)}\leqslant x\}q^{(r)}_{ij},\quad x\in\Omega,i=1,2,...,J^{(r)},r\in\mathbb{N}. \label{eqn:def omega}
\end{align}
Recall the definitions of operator $\mathcal{K}^{(r)}$ and kernel $\kappa^{(r)}$ in Lemma \ref{Lemma:IE reform I}. 
We can write $\kappa^{(r)}$ as:
\begin{align*}
    \kappa^{(r)}(x,u) =  \sum_{i=1}^{J^{(r)}}[\mathbf{1}\{x\geqslant 0\}+\bar{\omega}_i^{(r)}(x)]\cdot \mathbf{1}\{u\in(c^{(r)}_{i-1},c^{(r)}_i]\}
    =  \sum_{i=1}^{J^{(r)}}{\omega}_i^{(r)}(x)\cdot \mathbf{1}\{u\in(c^{(r)}_{i-1},c^{(r)}_i]\},
\end{align*}
for all $x,u\in\Omega,r\in\mathbb{N}$. 
Then we can also write operator $\mathcal{K}^{(r)}$ into the following equivalent form:
\begin{align}
    \mathcal{K}^{(r)}f(x) &=\int_\Omega \kappa^{(r)}(x,u)df(u)
    =\sum_{i=1}^{J^{(r)}}{\omega}_i^{(r)}(x)\cdot [f(c^{(r)}_i)-f(c^{(r)}_{i-1})], \quad \forall x\in\Omega, r\in\mathbb{N}. \label{eqn:K_r}
\end{align}
%Therefore, $\mathcal{K}^{(r)}$ is a finite-rank operator on $\bar{\mathbf{X}}$. 
According to Theorem {{}}5.2 of \cite{li2021numerical}, for a (finite-rank) operator $\mathcal{K}^{(r)}$ in the form of \eqref{eqn:K_r}, 
$(\mathcal{I}-\mathcal{K}^{(r)})^{-1}$ exists iff the matrix $H$ is non-singular,  where $H$ is defined by $H_{ij}:= \delta_{ij}-\omega^{(r)}_i(c^{(r)}_j)+\omega^{(r)}_{i+1}(c^{(r)}_j), \, i,j=1,2,...,J^{(r)}.$ Here $\delta_{ij}:=\mathbf{1}\{i=j\}\,(i,j\in\mathbb{N})$ and $\omega^{(r)}_{J^{(r)}+1}(x):=0,\forall x\in\mathbb{R}$. Then we only need to show $H$ is non-singular. By the definitions of $\{\omega^{(r)}_i\}_{i=1}^{J^{(r)}}$ in \eqref{eqn:def omega}, we have 
\begin{align*}
     H_{ij}= \begin{cases} \delta_{ij}-\sum_{k=1}^j q^{(r)}_{ik} +\sum_{k=1}^j q^{(r)}_{i+1,k} , \quad i,j\in \{1,2,...,J^{(r)}\}; i\neq J^{(r)},\\ \delta_{ij}-\sum_{k=1}^j q^{(r)}_{ik}-1, \quad i,j\in \{1,2,...,J^{(r)}\}; i= J^{(r)}.\end{cases} 
\end{align*}
Let us define %square matrix $E:E_{ij}= \mathbf{1}\{i\geqslant j\},\, i,j=1,2,...,J^{(r)}$ and  
another matrix $M$,  obtained by replacing the first row of $I-Q^{{(r)}\text{T}}$ with $\mathbf{1}^{\text{T}}\in \mathbb{R}^{J^{(r)}}$: 
\begin{align*}
    &M_{ij} = \begin{cases}\delta_{ij}- q^{(r)}_{ji},\quad i,j\in \{1,2,...,J^{(r)}\}; i \neq 1,  \\ 1,\quad i,j\in \{1,2,...,J^{(r)}\}; i = 1. \end{cases}
\end{align*}
Because the approximate finite-state MC is irreducible, rank$(M)=$rank$(I-Q^{{(r)}\text{T}})+1=J^{(r)}$. It is easy to see that $H^{\text{T}}$ can be transformed into $M$ by matrices' elementary row and column operations. Thus, $H$ is non-singular, $(\mathcal{I}-\mathcal{K}^{(r)})^{-1}$ exists and error bound \eqref{eqn:error bound naive} holds.  \Halmos
\end{proof}

According to Lemma \ref{Lemma:IE reform I}, error factor $||(\mathcal{I}-\mathcal{K}^{(r)})^{-1}||_O$ in \eqref{eqn:error bound naive} measures the approximate stationary distribution's sensitivity (in the inhomogeneous term $\xi$), while error factor $||\mathcal{K}p- \mathcal{K}^{(r)}p ||_\infty$ measures the kernel difference between the original and  approximate MCs.  

\subsection{Error Bound Computations.}
\label{sec:compuate error naive}
We complete the proof of Theorem \ref{Theorem: error bound naive} by showing that the operator-form error bound in Lemma \ref{Lemma:operator form error naive} is further bounded by $e_1\cdot e_2$, which is computable. 
The proof utilizes the following lemma.
\begin{lemma}[Theorem {{}}5.3, \cite{li2021numerical}]\label{Lemma:first factor naive}
For an instance of $\mathcal{K}^{(r)}$  in Lemma \ref{Lemma:operator form error naive}, 
$||(\mathcal{I}-\mathcal{K}^{(r)})^{-1}||_O^{-1}=\min_{k=0,1,...,J^{(r)}}\{y^*_{k}\}$, where $y^*_{k}$ is obtained from a linear program with $a_0\equiv 0$ $(k=0,1,...,J^{(r)})$:
\begin{align}
    &\min_{y_{k}\in\mathbb{R},\{a_j\}_{j=1}^{J^{(r)}} \in[-1,1]^{J^{(r)}} }     \quad \, y_{k}  \nonumber\\
    &s.t.  \quad\,\,  
     -y_{k}  \leqslant a_j+ \eta \delta_{jk}(1-a_j) -\sum_{i=1}^{J^{(r)}} \omega^{(r)}_i(c^{(r)}_j) [a_i-a_{i-1}] \leqslant y_{k},\, \quad j=0,1,...,J^{(r)}, \eta\in \{0,1\}.  \nonumber
\end{align}
\end{lemma}

\begin{proof}{Proof of Theorem \ref{Theorem: error bound naive}:}
We first show that in the operator-form error bound \eqref{eqn:error bound naive}, factor $||\mathcal{K}p- \mathcal{K}^{(r)}p ||_\infty\leqslant e_1$. 
Because $\left|\mathcal{K}{p}(x)- \mathcal{K}^{(r)}{p}(x) \right| = \left|\int_\Omega \kappa(x,u) dp(u) - \int_\Omega \kappa^{(r)}(x,u) dp(u) \right|    \leqslant  \int_\Omega \left| \kappa(x,u) -\kappa^{(r)}(x,u)\right| dp(u) \leqslant    \sup_{x,u\in\Omega} |\kappa(x,u) -\kappa^{(r)}(x,u)|,\,\forall x\in\Omega$, we have that 
\begin{align}
    ||\mathcal{K}{p}- \mathcal{K}^{(r)}{p} ||_\infty \leqslant \sup_{x,u\in\Omega} |\kappa(x,u) -\kappa^{(r)}(x,u)|\leqslant \sup_{x,u\in\Omega} |\bar{\kappa}(x,u) -\bar{\kappa}^{(r)}(x,u)|\leqslant e_1. \label{eq:second factor naive}
\end{align}
With regards to  factor $||(\mathcal{I}-\mathcal{K}^{(r)})^{-1}||_O$, it can be obtained by solving a series of linear programs as in Lemma \ref{Lemma:first factor naive}.
Because $\omega_i^{(r)}=\bar{\omega}^{(r)}_i+1$ (on domain $\Omega$) for $i=1,2,...,J^{(r)}$ according to \eqref{eqn:def omega}, we have that  $||(\mathcal{I}-\mathcal{K}^{(r)})^{-1}||_O= e_2$.  Therefore, $\left|\left| p-  p^{(r)} \right|\right|_\infty \leqslant ||(\mathcal{I}-\mathcal{K}^{(r)})^{-1}||_O \cdot ||\mathcal{K}p- \mathcal{K}^{(r)}p ||_\infty\leqslant e_1\cdot e_2.$
\Halmos
\end{proof}

%When we compute $e_2$, due to the structural similarities across $j$ and $k$,  methods including column generation can be applied to save computation time. 

\noindent\textbf{Remark 1 (weaker conditions for Theorem \ref{Theorem: error bound naive}).}\quad 
The proof of Lemma \ref{Lemma:operator form error naive} indicates that %the approximate finite-state MC's irreducibility is a sufficient condition for the error bound of \eqref{eqn:error bound naive}, while 
a necessary and sufficient condition for the error bound in \eqref{eqn:error bound naive} (or Theorem \ref{Theorem: error bound naive}) is the non-singularity of the matrix $M$ (rather than the approximate MC's irreducibility). $M$ is obtained by replacing the first row of $I-Q^{{(r)}\text{T}}$ with $\mathbf{1}^{\text{T}}\in \mathbb{R}^{J^{(r)}}$. 
This further implies the following necessary and sufficient condition for Theorem \ref{Theorem: error bound naive}: the approximate MC  has an absorbing communicating class. For a brief proof, under this condition, {\rm Ker}$(I-Q^{{(r)}\text{T}})=\{\alpha \pi^{(r)}|\alpha \in\mathbb{R}\}$, where $\pi^{(r)}$ is the stationary probability vector. Thus, $M$ is non-singular. On the other hand, if the approximate MC has two communicating classes not connected, then rank$(M)\leqslant$rank$(I-Q^{{(r)}\text{T}})+1\leqslant (J^{(r)}-2)+1 < J^{(r)}$, i.e., $M$ is singular.

\subsection{Errors in Steady-State Performance Measures.}\label{sec:performance error}
We now provide error bounds of performance measures computed by the approximate stationary distributions.
Such an error can be expressed as $\int_\Omega g(u) d\Delta p(u)=\int_\Omega g(u) d(p-p')(u)$, where $g$ is the performance function and $p'$ is the approximate stationary distribution.
Our error bound is based on the following observation: when we compute $\int_\Omega g(u) d\Delta p(u)$, the integral is bounded by the total variation $V(\Delta p)$ times the bound $||g||_\infty$. 
Using integral by parts, we switch the order of $\Delta p$ and $g$. 
In this way, the error is expected to be bounded by $V(g)$ times $||\Delta p||_\infty$: %. Formally, we have% the following result.
\begin{theorem}[Steady-State Performance Measure Errors]\label{theorem:performance measure error}
For distributions $p, p'\in\mathbf{D}$, and a measurable function $g(x)$ on $(\Omega, \mathcal{B})$ with a bounded total variation $V(g)<\infty$, we have that% performance measure error is bounded by
\begin{align}\label{ieq:performance measure error}
    \left|\int_\Omega g(x) dp(x) -\int_\Omega g(x) dp'(x)\right|\leqslant a V(g)\cdot ||p-p'||_\infty ,
\end{align}
where $a=2$.   If $g$ is continuous, $a=1$. 
\end{theorem}
\begin{proof}{Proof:} Without loss of generality, we assume $g$ is ``centralized'' such that $g(1)=0$. We next extend $g(x)$ such that $g(x)=g(0),\forall x< 0$ and $g(x)=g(1),\forall x>1$. 
Note that integrations in \eqref{ieq:performance measure error} are also valid in the sense of Lebesgue–Stieltjes measure: let $\hat{\mathrm{d}} h$ be the Lebesgue–Stieltjes measure associated with  $h:\mathbb{R}\mapsto \mathbb{R}$. Then  $
    \int_\Omega g(x)dp(x)  =  \int_\Omega g(x)\hat{\mathrm{d}} p(x)$ and $
    \int_\Omega g(x)dp'(x)  =  \int_\Omega g(x)\hat{\mathrm{d}} p'(x)$.
%The same result applies to $f'$. %Now we apply integral by parts under Lebesgue–Stieltjes measures. 

Define ${g}_C(x) = \frac{g(x{-})+g(x{+})}{2},\, g_R(x) = g(x+),\, {p}_C(x) = \frac{p(x{-})+p(x{+})}{2}$ and ${p'}_C(x) = \frac{p'(x{-})+p'(x{+})}{2}$
for all $x\in\Omega$. Then due to integral by parts (Theorem 21.67 in \cite{hewitt2013real}),  
$\int_\Omega {g}_C(x) \hat{\mathrm{d}}    p(x) +  \int_\Omega {p}_C(x) \hat{\mathrm{d}}  g_R(x)  =  {g}_C(1+){p}_C(1+)  -  {g}_C(0-){p}_C(0-)=0.$ Similarly, $\int_\Omega {g}_C(x) \hat{\mathrm{d}}    p'(x) +  \int_\Omega {p}'_C(x) \hat{\mathrm{d}}  g_R(x)=0.$ Thus, the steady-state performance measure error is bounded by 
\begin{align*}
    &\Big|\int_\Omega g(x) dp(x) -\int_\Omega g(x) dp'(x)\Big|\\ 
    =   &\Big| \int_\Omega \big[g(x)-{g}_C(x)\big] \hat{\mathrm{d}}    p(x) -  \int_\Omega {p}_C(x) \hat{\mathrm{d}}  g_R(x) - \int_\Omega \big[g(x)-{g}_C(x)\big] \hat{\mathrm{d}}    p'(x) + \int_\Omega {p'}_C(x) \hat{\mathrm{d}}  g_R(x)  \Big|\\
    \leqslant & \Big| \int_\Omega \big[g(x)-{g}_C(x)\big]  \hat{\mathrm{d}}  \big[p(x)-p'(x)\big] \Big| + \Big| \int_\Omega \big[{p}_C(x) -{p'}_C(x)\big]\hat{\mathrm{d}}  g_R(x) \Big|
\end{align*}
The second term is bounded by $V(g)\cdot ||p-p'||_\infty$. The first term is $0$ if $g$ is continuous. Otherwise, 
\begin{align*}
    &\Big| \int_\Omega \big[g(x)-{g}_C(x)\big]  \hat{\mathrm{d}}  \big[p(x)-p'(x)\big] \Big| = \Big| \int_{J_1\cup J_2} \big[g(x)-{g}_C(x)\big]  \hat{\mathrm{d}}  \big[p(x)-p'(x)\big] \Big|\\
    = & \Big| \sum_{x\in J_1\cup J_2} \big[g(x)-{g}_C(x)\big] \cdot \big[(p(x)-p(x-))-(p'(x)-p'(x-))\big]\Big|\\ \leqslant & \sum_{x\in J_1\cup J_2} \left|g(x)-{g}_C(x)\right| \cdot 2||p-p'||_\infty \leqslant \frac{V(g)}{2} \cdot 2||p-p'||_\infty \leqslant V(g) \cdot ||p-p'||_\infty,
\end{align*}
where $J_1= \{x\in \Omega\,|\,g(x{-})\neq g(x{+})\}$ and $J_2=\{x\in \Omega\,|\,g(x{-}) = g(x{+})\neq g(x) \}$. The transform from integral to summation in the second row is because $V(g)<\infty$, and sets $J_1$ and $J_2$ are both countable.
The inequality $\sum_{x\in J_1\cup J_2} \left|g(x)-{g}_C(x)\right|\leqslant \frac{V(g)}{2}$ is because every unit of difference in $|g(x)-{g}_C(x)|$ yields at least two units of total variation in $V(g)$. 
\Halmos
\end{proof}

\section{Consistency of Finite Approximation.}\label{sec:convergence naive}
%We consider the convergence of the approximate finite-state MCs' stationary distributions to that of the original MC. 
Theorem \ref{Theorem: error bound naive} has decomposed the finite approximation error into two factors: the kernel difference (original vs. approximate) and the approximate stationary distribution's sensitivity. 
In this section, we show that under certain regularity conditions, the sensitivity measure is bounded while the kernel difference measure converges to $0$. 
(Section \ref{sec:verify RC} will provide verifiable conditions that ensuring these regularity conditions.)
\begin{theorem}[Consistency in Stationary Distributions]\label{Theorem: convergence naive}
In Approximation \ref{A:Stationary Dist Eq}, if the original MC's kernel $\bar{\kappa}$ satisfies Condition \ref{A_UBV}-\ref{A_sol unique} given below, and the approximate MCs satisfy Condition \ref{A_T2} given below, then there exists $r_0\in\mathbb{N}$ such that for all $r>r_0$, stationary distribution $p^{(r)}$ is unique, the approximation error bound in \eqref{eqn:error bound naive} holds and the bound converges to $0$, i.e., $p^{(r)}\rightrightarrows p$ as $r\rightarrow \infty$.
\end{theorem}
The definitions and interpretations of Condition \ref{A_UBV}-\ref{A_T2} are as follows.
%These are referred to as ``regularity conditions" due to strong necessity and sufficiency (except Condition \ref{A_T2}, whose necessary and sufficient version can be found in \cite{li2021numerical} but is implicit, and Condition \ref{A_sol unique}).

Condition \ref{A_UBV}-\ref{A_MV} for the original MC's kernel $\bar{\kappa}$ are designed such that it is possible to find a finite set of ``representative" states to replace $\Omega$ when computing stationary distributions. 
Condition \ref{A_UBV} has been defined in Subsection \ref{sec:banach space}. 
The intuition behind Condition \ref{A_UBV} is that to make finite approximation possible, the original MC's transitions from the current state $u$ should vary within limit as $u$ varies from $0$ to $1$. %Only then the selected states can be ``representative" for the whole support.
Similarly, the original MC's transitions to the next state $x$ should vary within limit as $x$ varies from $0$ to $1$: this motivates the following regularity condition for kernel $\bar{\kappa}$. 
\begin{condition}[C\`adl\`ag and Countable Jump Discontinuities] \label{A_MV} Consider kernel $\tau\in{T}(\Omega, \mathcal{B})$. For all
$\epsilon>0$, there exists a finite split of $\Omega$ denoted by knots $0=s_1<s_2<...<s_{N_{\epsilon}}=1$ and intervals $E_1=[s_1,s_2),E_2=[s_2,s_3),...,E_{N_\epsilon}=\{s_{N_\epsilon}\}$ such that $\forall x_1,x_2\in E_i, i=1,2,...,N_\epsilon$, we have $\zeta_\tau(x_1,x_2)<\epsilon$ and $ v_\tau(x_1,x_2)<\epsilon$, where $\zeta_\tau(\cdot,\cdot)$ and $v_\tau(\cdot,\cdot)$ are distance functions respectively defined as $v_\tau(x_1,x_2):= V_u(\tau(x_2,u) - \tau(x_1,u))$ and $\zeta_\tau(x_1,x_2) := |\tau(x_2,1)-\tau(x_1,1)|,\forall x_1,x_2\in\Omega$. 
\end{condition}
For $\bar{\kappa}$ satisfying Condition \ref{A_MV}, the split knots divide support $\Omega$ into intervals, on each of which the original MC's transitions have a small variation as desired. 
Thus, for all $\epsilon$, the corresponding split knots defined by Condition \ref{A_MV} can be used as a group of ``representative" states. 
(Condition~\ref{A_MV} is referred to as ``c\`adl\`ag and countable jump discontinuities'' because it describes the property of limit existence, right continuity and countable jumps for $\tau(x,\cdot)$ respectively under the distance functions $v_\tau(\cdot,\cdot)$ and $\zeta_\tau(\cdot,\cdot)$ w.r.t. $x$ \citep{li2021numerical}.) 
%Moreover, the requirement for $\zeta_\tau(\cdot,\cdot)$ in Condition \ref{A_MV} is equivalent to $\tau(x,1)$ being right-continuous w.r.t. $x\in\Omega$. 

We also require that the original MC has a unique stationary distribution.
\begin{condition}[Unique Stationary Distribution]\label{A_sol unique} In Approximation \ref{A:Stationary Dist Eq}, the original MC's stationary distribution $p$ is unique and the solution set of \eqref{eqn:MC_BE_original} in $\mathbf{\bar{X}}$ is $\{\alpha p\}_{\alpha\in\mathbb{R}}$.
\end{condition}
%Verifiable sufficient conditions for Condition \ref{A_sol unique} are outlined in Section \ref{sec:verify RC}.

With regard to the approximate finite-state MC sequence, we assume that $(a)$ their states are ``representative" states as the split knots outlined in Condition \ref{A_MV}, %so that state $c^{(r)}_i (i=1,2,...,J^{(r)})$ will represent the original MC's all states in $[c^{(r)}_i,c^{(r)}_{i+1})\cap \Omega$,
and $(b)$ their transition matrices are defined by properly truncating the original MC's kernel. 
\begin{condition}[Proper Truncation]\label{A_T2}
In Approximation \ref{A:Stationary Dist Eq}, we have that:
%$\{\kappa^{(r)}\}_{r=1}^\infty$ in \eqref{Approx_kernel} specified by $\{\omega_i^{(r)}\}_{i=1}^{J^{(r)}}$, $\{c_i^{(r)}\}_{i=1}^{J^{(r)}} \,(r\in\mathbb{N})$. We have:
\begin{itemize}
    %\item[$(i)$] $\{\tau^{(r)} | r\in\mathbb{N}\}$ are defined as ``filtration finite approximation" and ``northwestern truncation and augmentation'' of $\tau$ in the sense that:
    \item[$(i)$](Representative states) There exist $\{\epsilon^{(r)}\}_{r\in\mathbb{N}}\subseteq \mathbb{R}_{+}$ and $\lim_{r\rightarrow \infty} \epsilon^{(r)}=0$ such that $\forall r \in \mathbb{N}$, $\forall x_1,x_2\in [c_{i}^{(r)},c_{i+1}^{(r)})\cap\Omega$, $i=1,2,...,J^{(r)}$, we have $ \zeta_{\bar{\kappa}}(x_1,x_2)<\epsilon^{(r)}$ and $ v_{\bar{\kappa}}(x_1,x_2)<\epsilon^{(r)}$. In other words, $\{c_i^{(r)}\}_{i=1}^{J^{(r)}}(r\in\mathbb{N})$ are defined by a sequence of finite splits for %the properties of c\`adl\`ag and countable jump discontinuities  of 
    $\bar{\kappa}$ in Condition \ref{A_MV} as $\epsilon\rightarrow 0$. 
    \item[$(ii)$] (Increasing partitions) $\{c_i^{(r_1)}\}_{i=1}^{J^{(r_1)}}\subseteq \{c_i^{(r_2)}\}_{i=1}^{J^{(r_2)}}$, $\forall r_1<r_2, r_1,r_2\in \mathbb{N}$. 
    \item[$(iii)$] (Kernel truncation) $q^{(r)}_{ij}=\bar{\kappa}(c^{(r)}_j,c^{(r)}_i)- \bar{\kappa}(c^{(r)}_{j-1},c^{(r)}_i)$,  $\forall i,j=1,2,...,J^{(r)}$, $r\in\mathbb{N}$.
    \item[$(iv)$] (Pointwise convergence) $
        \lim_{r\rightarrow\infty} \bar{\kappa}^{(r)}(x,u) = \bar{\kappa}(x,u),\forall u,x\in\Omega.$
\end{itemize}
\end{condition}
\begin{comment}
        \begin{align*}
            \omega_i^{(r)}(x)=\begin{cases}
            \sum_{j=1}^{J^{(r)}} \tau(c_{j-1}^{(r)},c_i^{(r)})\mathbf{1}\{x\in[c_{j-1}^{(r)},c_{j}^{(r)})\},\quad &x<c^{(r)}_{J^{(r)}}, \\
            \tau(+\infty, c_i^{(r)}),\quad &x\geqslant c^{(r)}_{J^{(r)}}.  
            \end{cases}
        \end{align*}
\end{comment}
%We refer to the definitions of approximate finite-state MC's states as  ``increasing partition'' because every approximate MC's states are included in that of subsequent MC. 
%Techniques for verifying Condition \ref{A_UBV}-\ref{A_T2} will be discussed in Section \ref{sec:verify RC}. 

Theorem \ref{theorem:performance measure error} and \ref{Theorem: convergence naive} imply the consistency of approximate steady-state performance measures.

\begin{corollary}[Consistency in Steady-State Performance Measures]
Given the settings and assumptions in Theorem \ref{Theorem: convergence naive} and a measurable performance function $g(x)$ on $(\Omega, \mathcal{B})$ such that $V(g)<\infty$, the steady-state performance measure associated with the approximate MC %'s stationary distribution $p^{(r)}$ 
converges to that associated with the original MC. %'s stationary distribution $p$. 
When we consider a set of performance functions  $\{g_i\}_{i\in S}$ on $(\Omega, \mathcal{B})$ such that $\sup_{i\in S}V(g_i)<\infty$, the steady-state performance measures associated with the approximate MC %$p^{(r)}$ 
uniformly converge to those associated with the original MC. %$p$. 
\end{corollary}

\subsection{Proof Outline.}
Our proof is based on the following lemma and balance equations' operator reformulation in Lemma \ref{Lemma:IE reform I}, which indicates
$f=p$ is a solution to $f=\xi+\mathcal{K}f$ while $f=p^{(r)}$ is a solution to $f=\xi+\mathcal{K}^{(r)}f$ $(r\in\mathbb{N})$.
\begin{lemma}[Lemma 2.3 of \cite{li2021numerical}]\label{lemma:anselone}
Consider function $h\in\bar{\mathbf{X}}$ and linear operators $\{\mathcal{K}\}\cup\{\mathcal{K}^{(r)}\}_{r\in\mathbb{N}}$ on $\bar{\mathbf{X}}$. Suppose $f=p$ and $f=p^{(r)}(r\in\mathbb{N})$ are respectively solutions to equations $f=h+\mathcal{K}f$ and $f=h+\mathcal{K}^{(r)}f$ $(r\in\mathbb{N})$ on $\bar{\mathbf{X}}$. If  $\{\mathcal{K}\}\cup\{\mathcal{K}^{(r)}\}_{r\in\mathbb{N}}$ further satisfy: 
\begin{itemize}
    \item[$(A1)$] $\{\mathcal{K}^{(r)}\}_{r\in\mathbb{N}}$ are \textit{collectively compact}, i.e., the set $\{\mathcal{K}^{(r)}f\,|\,r\in\mathbb{N}, f\in \mathbf{\bar X},||f||_\infty\leqslant 1 \}$ is relatively compact,
    \item[$(A2)$] $\{\mathcal{K}^{(r)}\}_{r\in\mathbb{N}}$ are \textit{consistent operators} to $\mathcal{K}$, i.e., $\mathcal{K}^{(r)}f \rightrightarrows \mathcal{K}f,\forall  f\in \mathbf{\bar X}$, and
    \item[$(A3)$] $(\mathcal{I}-\mathcal{K})^{-1}$ exists, 
\end{itemize}
then  we have the following results:
\begin{itemize}
    \item[$(B1)$] there exists $r_0\in\mathbb{N}$ such that $\forall r>r_0$, $p^{(r)}$ is unique and follows the error bound in \eqref{eqn:error bound naive},
    \item[$(B2)$] consistency, i.e., $p^{(r)} \rightrightarrows p$, and 
    \item[$(B3)$] stability, i.e.,  $\{||(\mathcal{I}-\mathcal{K}^{(r)})^{-1}||_O\,|\, r>r_0\}$ is bounded.
\end{itemize}
\end{lemma}
Note $(B1)$-$(B2)$ include all results of Theorem \ref{Theorem: convergence naive}. Thus, we only need to show $(A1)$-$(A3)$ to prove Theorem \ref{Theorem: convergence naive} according to Lemma \ref{lemma:anselone}.
Our proof can be completed in two steps: in Subsection \ref{sec:consistency and collective compactness}, we prove that Conditions \ref{A_UBV}, \ref{A_MV} and \ref{A_T2} imply collective compactness ($A1$) and operator consistency ($A2$). In Subsection \ref{sec:invertibility}, we prove that Condition \ref{A_sol unique} implies invertibility ($A3$). 

Particularly, we show that under appropriate assumptions,  our finite approximation approach is near-optimal in approximating the original MC's stationary distribution among all approximation methods using discrete distributions (Subsection \ref{sec:near optimality}). 

\subsection{Collectively Compact and Consistent Operators.}\label{sec:consistency and collective compactness}
We use  the (finite-rank) structure of $\{\mathcal{K}^{(r)}\}_{r\in\mathbb{N}}$ described in \eqref{eqn:K_r} and the following Lemma to show that Condition \ref{A_UBV},\ref{A_MV} and \ref{A_T2} imply collective compactness $(A1)$ and operator consistency $(A2)$.  
\begin{lemma}[Theorem 3.5 in \cite{li2021numerical}]\label{lemma:cc and c}
For operators $\{\mathcal{K}\}\cup\{\mathcal{K}^{(r)}\}_{r\in\mathbb{N}}$ as defined in Lemma \ref{Lemma:IE reform I},  $\{\mathcal{K}^{(r)}\}_{r\in\mathbb{N}}$ are collectively compact and consistent operators iff 
\begin{itemize}
    \item[$(i)$] (Pointwise convergence) $\forall f\in \bar{\mathbf{X}}$, $\forall  x\in\Omega$, $\lim_{k\rightarrow\infty} \mathcal{K}^{(r)}f(x) = \mathcal{K} f(x).$
    \item[$(ii)$] (Uniform c\`adl\`ag and countable jump discontinuities) for all $\epsilon>0$, there exists a finite split of $\Omega$ denoted by knots $0 =s_1<s_2<\ldots <s_{N_{\epsilon}}=1$ and intervals $E_1=[s_1,s_2),E_2=[s_2,s_3),\ldots, E_{N_\epsilon}=\{s_{N_\epsilon}\}$ such that $\forall r\in\mathbb{N},  x_1,x_2\in E_i,  i=1,2,\ldots ,N_\epsilon$, we have $\sum_{i=1}^{J^{(r)}} |\Delta \omega^{(r)}_i(x_2)-\Delta \omega^{(r)}_i(x_1)|\leqslant \epsilon$. Here $\Delta \omega^{(r)}_i(x) := \omega^{(r)}_{i+1}(x) - \omega^{(r)}_{i}(x), x\in\Omega,i=1,2,\ldots ,J^{(r)}, r\in\mathbb{N}$;  %$\{\omega_i^{(r)}\}_{i=1}^{J^{(r)}}(r\in\mathbb{N})$ are as defined in \eqref{eqn:def omega}; 
    and particularly, $\omega^{(r)}_{J^{(r)}+1}(x)\equiv 0, \forall r \in \mathbb{N}$.
\end{itemize}
\end{lemma}

\begin{proof}{Proof of (A1)-(A2) in Lemma \ref{lemma:anselone}:} %By Lemma \ref{lemma:cc and c}, we only need to prove $\{\mathcal{K}^{(r)}\}_{r\in\mathbb{N}}$ satisfy pointwise convergence and uniform c\`adl\`ag and countable jump discontinuities. 

(\textit{Pointwise convergence.})\quad  We use dominated convergence. 
Recall definitions of $\kappa=\bar{\kappa}+1$ and $\kappa^{(r)}=\bar{\kappa}^{(r)}+1 (r\in\mathbb{N})$ (on domain $\Omega^2$) in Lemma \ref{Lemma:IE reform I}. 
Then $(iv)$ of Condition \ref{A_T2} implies $\lim_{r\rightarrow \infty} \kappa^{(r)}(x,u)=\kappa(x,u),\forall x,u\in\Omega$. 
For all $f\in \mathbf{X}$, $\mathcal{K}^{(r)}f$ can also be written as $ \mathcal{K}^{(r)} F (x)=  \int_\Omega \kappa^{(r)} (x,u)df(u), \, \forall x\in\Omega, r\in\mathbb{N}.$
Thus, for all $f\in\mathbf{D}$ and $x\in\Omega$, by the dominated convergence, we have $\lim_{r\rightarrow \infty} \mathcal{K}^{(r)}f(x) =\int_\Omega \lim_{r\rightarrow \infty} \kappa^{(r)} (x,u)df(u) =  \int_\Omega \kappa (x,u)df(u) =  \mathcal{K}f(x).$ 
By the linearity, for all $f\in \mathbf{X}$, $\lim_{r\rightarrow \infty} \mathcal{K}^{(r)}f(x) =   \mathcal{K}f(x), \forall x \in \Omega$. 

We next expand this result to $f\in\bar{\mathbf{X}}$.
Consider $f=f'+e$, $f\in\bar{\mathbf{X}}, f'\in \mathbf{X}$ and $||e||_\infty<\epsilon$. 
Here $f'$ can be regarded as an approximation to $f$ in the subspace of $\mathbf{X}$ with a bounded error $e$. 
For all $x\in\Omega$, we have already proven $\lim_{r\rightarrow \infty} \mathcal{K}^{(r)}f'(x) = \mathcal{K}f'(x)$. 
We only need to prove $\lim_{\epsilon\rightarrow 0}\limsup_{r\rightarrow \infty} |\mathcal{K}^{(r)}e(x) - \mathcal{K}e(x)|=0$. 
Assume there exists a constant $C<\infty$, whose existence will be proven later, such that $||\mathcal{K}^{(r)}h - \mathcal{K}h||_\infty<C||h||_\infty$ for all $h\in\bar{\mathbf{X}},r\in\mathbb{N}$. 
Then $\lim_{\epsilon\rightarrow 0}\limsup_{r\rightarrow \infty} |\mathcal{K}^{(r)}e(x) - \mathcal{K}e(x)|\leqslant \lim_{\epsilon\rightarrow 0} C\epsilon =0$. 
Therefore, $\lim_{r\rightarrow \infty} \mathcal{K}^{(r)}f(x) = \mathcal{K}f(x)$.

Lastly, we prove that there exists a constant $C<\infty$ such that $||\mathcal{K}^{(r)}h - \mathcal{K}h||_\infty<C||h||_\infty$ for all $h\in\bar{\mathbf{X}},r\in\mathbb{N}$. According to equation (SM2.8) in Theorem 3.1 of \cite{li2021numerical}, for any kernel $\tau\in T(\Omega,\mathcal{B})$ satisfying Condition \ref{A_UBV} and its associated operator $\mathcal{L}$ in Definition \ref{Def:K}, we have $||\mathcal{L}h||_\infty \leqslant ||\mathcal{L}||_O||h||_\infty \leqslant ||h||_\infty\cdot \sup_{x\in\Omega} [2V_u(\tau(x,u))+ \tau(x,1)]$ for all $h\in\bar{\mathbf{X}}$. 
Recall operators $\mathcal{K}$ and $\mathcal{K}^{(r)}(r\in\mathbb{N})$ are respectively associated with kernels $\kappa$ and $\kappa^{(r)}(r\in\mathbb{N})$. Then we have
\begin{align*}
    &||\mathcal{K}h-\mathcal{K}^{(r)}h||_\infty\leqslant ||\mathcal{K}h||_\infty+||\mathcal{K}^{(r)}h||_\infty\\
    \leqslant & ||h||_\infty\cdot  \sup_{x\in \Omega} \Big[  2V_u(\kappa(x,u))+ \kappa(x,1) \Big]+||h||_\infty\cdot  \sup_{x\in \Omega} \Big[  2V_u(\kappa^{(r)}(x,u))+ \kappa^{(r)}(x,1) \Big]\\ \leqslant & ||h||_\infty\cdot  \sup_{x\in \Omega} \Big[  2V_u(\kappa(x,u))+ 2 \Big]+||h||_\infty\cdot  \sup_{x\in \Omega} \Big[  2V_u(\kappa(x,u))+ 2 \Big].   
\end{align*}
Note here we use the facts $\kappa(x,u)\leqslant \bar{\kappa}(x,u)+1 \leqslant 2, \kappa^{(r)}(x,u)\leqslant \bar{\kappa}^{(r)}(x,u)+1\leqslant 2, \forall x,u\in \Omega$; and $\sup_{x\in \Omega} V_u(\kappa^{(r)}(x,u))=\sup_{x\in \Omega}V_u(\bar{\kappa}^{(r)}(x,u)) \leqslant \sup_{x\in \Omega}V_u(\bar{\kappa}(x,u))=\sup_{x\in \Omega}V_u(\kappa(x,u))$ due to $(iii)$  in Condition \ref{A_T2}. Thus, we can define $C:=4\sup_{x\in \Omega} V_u(\kappa(x,u))+4=4\sup_{x\in \Omega} V_u(\bar\kappa(x,u))+4$. Because $\bar{\kappa}$ satisfies Condition \ref{A_UBV}, $C<\infty$. 

(\textit{Uniform c\`adl\`ag and countable jump discontinuities.})\quad  
Because $\kappa=\bar{\kappa}+1$ (on domain $\Omega^2$), the kernel $\bar{\kappa}$ in $(i)$ of Condition \ref{A_T2} can be replaced by $\kappa$. 
Recall $\epsilon^{(r)}(r\in\mathbb{N})$ is associated with $\{c_i^{(r)}\}_{i=1}^{J^{(r)}}$ in $(i)$ of Condition \ref{A_T2}. 
For all $\epsilon>0$, we show that with $r_1 = \inf \{r\in\mathbb{N}|\epsilon^{(r)}< \frac{\epsilon}{2}\}$, knots $\{c_i^{(r_1)}\}_{i=1}^{J^{(r_1)}}$ form a finite split satisfying $(ii)$ of Lemma \ref{lemma:cc and c}  with $\epsilon$. 
Indeed, $(a)$ for $r\leqslant r_1$, $\forall x_1,x_2\in [c_{i}^{(r_1)},c_{i+1}^{(r_1)})\cap\Omega,$ $i=1,2,...,J^{(r_1)}$, we have that $x_1,x_2\in[c_{i}^{(r)},c_{i+1}^{(r)})\cap\Omega$ for some $i\in\{1,2,...,J^{(r)}\}$ (due to the increasing partition property in $(ii)$ of Condition \ref{A_T2}) and $\omega_i^{(r)}(x_2)=\omega_i^{(r)}(x_1)\,(i=1,2,...,J^{(r)}+1)$. 
Thus, $\sum_{i=1}^{J^{(r)}} \big|\Delta\omega_i^{(r)}(x_2)-\Delta\omega_i^{(r)}(x_1) \big|=0.$
(b) For $r> r_1$, we first define function 
$c^{(r)}(x) :=\sum_{i=1}^{J^{(r)}} c_i^{(r)}\mathbf{1}\{x\in[c^{(r)}_i,c^{(r)}_{i+1})\}, x\in\Omega.$
Then, $\forall x_1,x_2\in [c_{i}^{(r_1)},c_{i+1}^{(r_1)})\cap\Omega,$ $i=1,2,...,J^{(r_1)}$, we have that
\begin{align*}
    &\sum_{i=1}^{J^{(r)}} \big|\Delta\omega_i^{(r)}(x_2)-\Delta\omega_i^{(r)}(x_1) \big|
    = V_u\big( \kappa^{(r)}(x_2,u)-\kappa^{(r)}(x_1,u) \big) + \big|\kappa^{(r)}(x_2,1)-\kappa^{(r)}(x_1,1)\big| \\
    = & V_u\big( \kappa^{(r)}(c^{(r)}(x_2),u)-\kappa^{(r)}(c^{(r)}(x_1),u) \big) + \big|\kappa^{(r)}(c^{(r)}(x_2),1)-\kappa^{(r)}(c^{(r)}(x_1),1)\big| \\
    %= & V_u\big( \kappa^{(r)} (c^{(r)}(x_2),u)-\kappa^{(r)}(c^{(r)}(x_1),u) \big)\\
    \leqslant & V_u\big( \kappa (c^{(r)}(x_2),u)-\kappa(c^{(r)}(x_1),u) \big) +\big|\kappa(c^{(r)}(x_2),1)-\kappa(c^{(r)}(x_1),1)\big| \quad (\text{due to $(iii)$ of Condition \ref{A_T2}} )\\
    \leqslant & v_\kappa (c^{(r)}(x_2),c^{(r)}(x_1))+ \zeta_\kappa(c^{(r)}(x_2),c^{(r)}(x_1))
    \leqslant  2 \epsilon^{(r_1)} \leqslant   \epsilon,
\end{align*}
where the last row follows $c^{(r)}(x_2),c^{(r)}(x_1)\in[c_{i}^{(r_1)},c_{i+1}^{(r_1)})\cap \Omega$ for some $i\in\{1,2,...,J^{(r_1)}\}$ due to the increasing partition property in $(ii)$ of Condition \ref{A_T2}. Recall  $\forall x_1,x_2\in [c_{i}^{(r_1)},c_{i+1}^{(r_1)})\cap\Omega$, $i=1,2,...,J^{(r_1)}$, we have $v_\kappa(x_1,x_2)+\zeta_\kappa(x_1,x_2)\leqslant 2\epsilon^{(r_1)}$. Thus, $v_\kappa (c^{(r)}(x_2),c^{(r)}(x_1))+ \zeta_\kappa(c^{(r)}(x_2),c^{(r)}(x_1)) \leqslant   2\epsilon^{(r_1)} \leqslant   \epsilon$. 
Lastly, due to $(a)$ and $(b)$, the finite split $\{c_i^{(r_1)}\}_{i=1}^{J^{(r_1)}}$ satisfies $(ii)$ of Lemma \ref{lemma:cc and c} with $\epsilon$.
\Halmos
\end{proof}

\subsection{Invertibility.}\label{sec:invertibility}
We use the following Lemma to show that Condition \ref{A_sol unique} implies $(A3)$. 

\begin{lemma}[Lemma 2.4 and Theorem 3.4 in \cite{li2021numerical}]\label{lemma:invertibility}
For kernel $\kappa$ and its associated  operator $\mathcal{K}$ in Lemma \ref{Lemma:IE reform I}, if $\kappa$ satisfies Condition \ref{A_UBV}-\ref{A_MV} and $\mathcal{K}$ satisfies ${\rm Ker}(\mathcal{I}-\mathcal{K})=\{0\}$, then $(\mathcal{I-K})^{-1}$ exists.   
\end{lemma}
This lemma is based on the Riesz-Schauder's theorem: Condition \ref{A_UBV}-\ref{A_MV} for $\kappa$ imply the compactness of $\mathcal{K}$, and ${\rm Ker}(\mathcal{I}-\mathcal{K})=\{0\}$ implies the existence of $(\mathcal{I-K})^{-1}$ by the Riesz-Schauder's theorem. %See a detailed proof in \cite{li2021numerical}. 

\begin{proof}{Proof of $(A3)$ in Lemma \ref{lemma:anselone}:}
%According to Lemma \ref{lemma:invertibility}, to prove $(A3)$, i.e., the existence of $(\mathcal{I-K})^{-1}$, 
We only need to show $(i)$ $\kappa$ satisfies Condition \ref{A_UBV}-\ref{A_MV} and $(ii)$ ${\rm Ker}(\mathcal{I}-\mathcal{K})=\{0\}$. Verifying $(i)$ is not hard since $\kappa=\bar{\kappa}+1$ (on domain $\Omega^2$) and $\bar{\kappa}$ satisfies Condition \ref{A_UBV}-\ref{A_MV}. 

To prove $(ii)$, we only need to show that \eqref{eqn:MC_OE_original} has a unique solution in $\bar{\mathbf{X}}$. To prove the solution uniqueness, we use the fact that \eqref{IE reform 1: right} is reversible, i.e., for all $f\in\bar{\mathbf{X}}$, we have
\begin{align}\label{IE reform 2: rightleft}
    f(1)=1; f(x)= \int_\Omega \bar{\kappa}(x,u)df(u),\forall x\in\Omega. \quad \Leftrightarrow \quad f(x)=\xi(x)+ \int_\Omega {\kappa}(x,u)df(u),\forall x\in\Omega. 
\end{align}
(The proof will be presented shortly.) The system of equations on the right hand side (RHS) is exactly \eqref{eqn:MC_IE_original}, while the left hand side (LHS) has a unique solution $p\in\bar{\mathbf{X}}$ due to Condition \ref{A_sol unique}. Thus, \eqref{eqn:MC_IE_original} as well as its operator form \eqref{eqn:MC_OE_original} has a unique solution $p \in \bar{\mathbf{X}}$. 

Lastly, we show the equivalence in \eqref{IE reform 2: rightleft} to complete the proof. 

($\Rightarrow$) Consider any fixed $f\in\bar{\mathbf{X}}$. 
The LHS of \eqref{IE reform 2: rightleft} implies $f(x)=\mathbf{1}\{x\geqslant 0\}\cdot[f(1)-1] + \int_\Omega \bar{\kappa}(x,u)df(u)= \xi(x)+ f(1) + \int_\Omega \bar{\kappa}(x,u)df(u),\forall x\in\Omega$. 
Because $\kappa=\bar{\kappa}+1$ (on domain $\Omega^2$) according to Lemma \ref{Lemma:IE reform I}, we only need to prove $f(1)=\int_\Omega df(u)$ to obtain the RHS. 
Indeed, if $f\in\mathbf{X}$, then $f(1)=\int_\Omega df(u)$ because $f$ defines a finite signed measure. 
If $f\in\bar{\mathbf{X}}\backslash \mathbf{X}$, then $f(1)=\lim_{k\rightarrow \infty}f_k(1)= \lim_{k\rightarrow\infty} \int_\Omega df_k(u)=\int_\Omega df(u)$, where $\{f_k\}_{k\in\mathbb{N}}\subseteq \mathbf{X}$ is any sequence converging to $f$ and the last equality is due to \eqref{Def:K_2}.

($\Leftarrow$) Consider any fixed $f\in\bar{\mathbf{X}}$. We have shown $f(1)=\int_\Omega df(u)$ in ($\Rightarrow$). Thus, the RHS implies 
\begin{align}
    f(x)= \xi(x)+ \int_\Omega df(u) + \int_\Omega \bar{\kappa}(x,u)df(u) = -1 + f(1) + \int_\Omega \bar{\kappa}(x,u)df(u),\quad \forall x\in\Omega.    \label{eqn:142415}
\end{align}
Consider $x=1$. We have $f(1)=-1 + f(1) + \int_\Omega df(u)= -1+2f(1)$. Thus, $f(1)=1$. Plugging this into \eqref{eqn:142415} yields $f(x) = \int_\Omega \bar{\kappa}(x,u)df(u),\forall x\in\Omega$, which constitutes the LHS with $f(1)=1$.  \Halmos 
\end{proof}

Now we complete the proof of Theorem \ref{Theorem: convergence naive} with results in Subsection \ref{sec:consistency and collective compactness}-\ref{sec:invertibility}.
\begin{proof}{Proof of Theorem \ref{Theorem: convergence naive}:}\quad 
In Subsection \ref{sec:consistency and collective compactness}-\ref{sec:invertibility}, we have shown $(A1)$-$(A3)$ in Lemma \ref{lemma:anselone} hold. Thus, $(B1)$-$(B2)$ in Lemma \ref{lemma:anselone} also hold, i.e., we have the results of Theorem \ref{Theorem: convergence naive}. \Halmos 
\end{proof}

\subsection{Near-Optimality of  Finite Approximation.}\label{sec:near optimality}
%\noindent\textbf{Remark 2: near-optimal approximation.} 
We make the following assumption.
%With the settings in Theorem \ref{Theorem: convergence naive} and the following extra assumption, our finite approximation approach is near-optimal in approximating the original MC's stationary distribution among all methods using discrete distributions.
\begin{assumption}\label{A_Lipschitz}
In Approximation \ref{A:Stationary Dist Eq}, we have
$(i)$  kernel $\bar{\kappa}$ is Lipschitz continuous on $(x,u)\in\Omega^2$, and $(ii)$ there exists a small interval $[x_1,x_2]\subseteq \Omega$ on which $p$ is continuous and strictly increasing.
\end{assumption}
We let $\delta:=p(x_2)-p(x_1) >0$ denote the steady-state probability associated with the small interval in Assumption \ref{A_Lipschitz}.
Without loss of generality, we assume the approximate MCs' states $\{c^{(r)}_i\,|\,i=1,2,...,J^{(r)},r\in\mathbb{N}\}$ are dense in $\Omega$ and $\max_{i=2,3,...,J^{(r)}}|c^{(r)}_i-c^{(r)}_{i-1}|=O(\frac{1}{J^{(r)}})$. 
If the denseness or the convergence rate does not hold, we may replace the definition of $\{c_i^{(r)}\}_{i=1}^{J^{(r)}}$ with $\{c_i^{(r)}\}_{i=1}^{J^{(r)}}\cup \{\frac{i}{2^{\sigma(r)}}\,|\,i=0,1,...,2^{\sigma(r)}\}$, where $\sigma(r)=\inf\{i\in\mathbb{N}\,|\,2^i+1\geqslant J^{(r)},i\geqslant r\}\,(r\in\mathbb{N})$. 
%, with which the approximate MC sequence still satisfies Condition \ref{A_T2} and the above denseness and convergence rate also hold.

The near-optimality of finite approximation is proven below:
\begin{theorem}[Near-optimality of Finite Approximation]
\label{theorem:near-optimality}
Consider the settings and assumptions in Theorem \ref{Theorem: convergence naive}, the  finite approximate distributions $\{p^{(r)}\}_{r\in\mathbb{N}}$, which respectively have ${J}^{(r)}$ jumps, and any other sequence of discrete distributions $\{\hat{p}^{(r)}\}_{r\in\mathbb{N}}\subseteq\mathbf{D}$, which respectively have $\hat{J}^{(r)}$ jumps. If Assumption \ref{A_Lipschitz} holds and  $\lim_{r\rightarrow\infty}\frac{\hat{J}^{(r)}}{{J}^{(r)}}=1$, i.e., the two sequences of discrete distributions  have asymptotically the same number of jumps, then ${||p^{(r)}-p||_\infty}=O({||\hat p^{(r)}-p||_\infty})$, i.e., as approximators to $p$, $\{p^{(r)}\}_{r\in\mathbb{N}}$ are at least as accurate as $\{\hat{p}^{(r)}\}_{r\in\mathbb{N}}$ by a constant multiplier.
\end{theorem}
\begin{proof}{Proof:}
We only need to show $|| p^{(r)}-p||_\infty=O(\frac{1}{{J}^{(r)}})$ while 
$||\hat p^{(r)}-p||_\infty\geqslant \frac{\delta}{2(\hat{J}^{(r)}+1)}$. 
Then ${||p^{(r)}-p||_\infty}=O({||\hat p^{(r)}-p||_\infty})$ due to $\lim_{r\rightarrow\infty}\frac{\hat{J}^{(r)}}{{J}^{(r)}}=1$. 

To bound $||\hat p^{(r)}-p||_\infty$, we only need to consider the small interval in Assumption \ref{A_Lipschitz}, on which $p$ is continuous and strictly increasing. For any discrete approximate distribution with $J$ jumps, %the maximal jump size on this small interval is at least $\frac{\delta}{J}$. %: we only need to consider the small interval in Assumption \ref{A_Lipschitz}. %Because $p$ is continuous , the approximation error on this interval is at least $\frac{\delta}{2J}$, and the overall error is at least $O(\frac{1}{J})$.
the approximation error on this interval is at least $\frac{\delta}{2(J+1)}$ because  $p$ is continuous on this interval. Then the overall error is also at least $\frac{\delta}{2(J+1)}$. Thus, $||\hat p^{(r)}-p||_\infty\geqslant \frac{\delta}{2(\hat{J}^{r}+1)}$. 

Lastly, we consider $|| p^{(r)}-p||_\infty$. 
%If we apply finite approximation with the settings of Theorem \ref{Theorem: convergence naive}, 
The approximation error bound is given in \eqref{eqn:error bound naive} of Lemma \ref{Lemma:operator form error naive}. 
Factor $||(\mathcal{I}-\mathcal{K}^{(r)})^{-1}||_O$ is uniformly bounded according to $(B3)$ of Lemma \ref{lemma:anselone}. 
Factor $||\mathcal{K}p- \mathcal{K}^{(r)}p ||_\infty$ is bounded by $\sup_{x,u\in\Omega} |\bar{\kappa}(x,u) -\bar{\kappa}^{(r)}(x,u)|$ according to \eqref{eq:second factor naive}. 
Recall the definition of $\bar{\kappa}^{(r)}$ under Condition \ref{A_T2}, which can be regarded as a step-wise function approximating $\bar{\kappa}$. 
Because $\bar{\kappa}(x,u)$ is Lipschitz continuous and  $\max_{i=2,3,...,J^{(r)}}|c^{(r)}_i-c^{(r)}_{i-1}|=O(\frac{1}{J^{(r)}})$, we have that $\sup_{x,u\in\Omega} |\bar{\kappa}(x,u) -\bar{\kappa}^{(r)}(x,u)|= O(\frac{1}{J^{(r)}})$. 
In this way, the approximation error $|| p^{(r)}-p||_\infty$ is also $O(\frac{1}{J^{(r)}})$. \Halmos
\end{proof}

According to Theorem \ref{theorem:near-optimality}, when we compute the original MC's stationary distribution, the finite approximate distributions are at least as accurate as the empirical distributions generated by a simulation approach (e.g., discrete event simulation, MCMC, etc.) by a constant multiplier because the later are also discrete distributions.

\section{Verifiable Sufficient Conditions.}\label{sec:verify RC}
%In this section, we provide sufficient conditions that ensure the abstract Condition \ref{A_UBV}-\ref{A_T2}.
\subsection{Sufficient Conditions for Conditions 1,2 and 4.}\label{sec:verify c1-3}
%Condition \ref{A_UBV} requires a uniformly bounded total variation of $\tau(x,u)$ along $u$ for every fixed $x$. 
Suppose $\tau(x,u)$ is Lipschitz continuous w.r.t. $u$, i.e., there exists $L>0$ such that $|\tau(x,u_2)-\tau(x,u_1)|\leqslant L|u_2-u_1|,\forall x,u_1,u_2\in\Omega$. Then Condition \ref{A_UBV} holds because $\forall x\in\Omega$, we have 
\begin{align}
    V_u(\tau(x,u))= \sup_{0\leqslant u_1<u_2<...<u_n\leqslant 1; n\in\mathbb{N} } \sum_{i=2}^n |\tau(x,u_i)-\tau(x,u_{i-1})|\leqslant \sup_{0\leqslant u_1<u_2<...<u_n\leqslant 1; n\in\mathbb{N} } \sum_{i=2}^n L|u_i-u_{i-1}|\leqslant L. \nonumber
\end{align}
A stronger sufficient condition for Condition \ref{A_UBV} is that the partial derivative $\frac{\partial\tau(x,u)}{\partial u}$ is bounded or continuous for $(x,u)\in\Omega^2$, which yields $V_u(\tau(x,u)) = \int_{\Omega} \left|\frac{\partial\tau(x,u)}{\partial u}\right| du \leqslant \sup_{x,u\in\Omega}\left|\frac{\partial\tau(x,u)}{\partial u}  \right|<\infty.$

%\subsection{Verify Condition 2 of c\`adl\`ag and countable jump discontinuities}
Condition \ref{A_MV} can be verified similarly. % via a similar approach to Condition \ref{A_UBV}. 
If $\left|\frac{\partial^2\tau(x,u)}{\partial x \partial u}\right|<L'$ for all $(x,u)\in\Omega^2$, then for all $x_1,x_2\in\Omega$, 
\begin{align}
    v_\tau(x_1,x_2) = V_u(\tau(x_2,u)-\tau(x_1,u)) = \int_\Omega \left|\frac{\partial[\tau(x_2,u)-\tau(x_1,u)]}{\partial u}\right| du = \int_\Omega L'|x_2-x_1| du  \leqslant L'|x_2-x_1|.\nonumber %\label{eqn:verify c2_1} 
\end{align}
If $|\frac{d\tau(x,1)}{d x}| \leqslant L''$ for all $x\in\Omega$, then $\zeta_\tau(x_1,x_2) = |\tau(x_2,1)-\tau(x_1,1)| \leqslant L''|x_2-x_1|$ for all $x_1,x_2\in\Omega$.  
Let $L'''=:\max\{L',L''\}$. Then Condition \ref{A_MV} is satisfied: for every $\epsilon>0$, we may define the corresponding split knots as $\big\{0,\frac{1}{2^n},\frac{2}{2^n},..., 1\big\}$, where $n= \inf\{i\in\mathbb{N}_{++} \,|\,\frac{L'''}{2^i}\leqslant \epsilon\}$. 

Regarding Condition \ref{A_T2}, $(i)$ simply requires that the approximate MCs' states are specified in accordance with the finite splits in Condition \ref{A_MV}. 
%, which is not hard to verify. 
For $(ii)$, if  $\{c_i^{(r_1)}\}_{i=1}^{J^{(r_1)}}$ is not a subset of $\{c_i^{(r_2)}\}_{i=1}^{J^{(r_2)}}$, we may replace  $\{c_i^{(r_2)}\}_{i=1}^{J^{(r_2)}}$ with $\{c_i^{(r_1)}\}_{i=1}^{J^{(r_1)}}\cup\{c_i^{(r_2)}\}_{i=1}^{J^{(r_2)}}$. 
Item $(iii)$ requires that the transition matrices of the approximate MCs are truncated from the original MC's transition kernel. 
The verification of the limit in  $(iv)$ is not hard by using properties of $\bar{\kappa}$ (e.g., the continuity of $\bar{\kappa}$ may imply $(iv)$; see an example in Subsection \ref{sec:GG1G}). 

\subsection{Sufficient Condition for Condition 3.}\label{sec:verify IE sol unique}
We first define a class of  accessible singletons.
\begin{definition} [Independent Mass Singleton]\label{def:independent mass singleton}
For an original MC in Approximation \ref{A:Stationary Dist Eq}, a singleton $\{x_0\}$ with state $x_0\in\Omega$ will be called an independent mass singleton if after replacing the definition of $\bar{\kappa}(\cdot,x_0)$ with an arbitrary $h\in\mathbf{D}$, the MC with the modified kernel satisfies: 
\begin{itemize}
    \item[$(i)$] a stationary distribution exists, and
    \item[$(ii)$] every stationary distribution includes a mass at $x_0$.
\end{itemize} 
\end{definition}
If $\{x_0\}$ is an independent mass singleton, then for all $h\in\mathbf{D}$, there exists $p'\in\mathbf{D}$ such that $\int_{\{x_0\}}dp'(u)>0$ and $p'(x) = \int_\Omega \bar{\kappa}' (x,u) dp'(u),\forall x\in\Omega$, where the modified kernel 
\begin{align}
    \bar\kappa'(x,u) := \begin{cases}\bar{\kappa}(x,u),\quad u\neq x_0, \\h(x), \quad u=x_0. \end{cases}  \label{eqn:kernel modified}
\end{align}

\begin{theorem}[Solution Uniqueness by Independent Mass Singleton]
\label{theorem: sol unique}
In Approximation \ref{A:Stationary Dist Eq}, if $(i)$ kernel $\bar{\kappa}$ satisfies Condition \ref{A_UBV}-\ref{A_MV} and $(ii)$ there exists an $x_0\in\Omega$ that forms an independent mass singleton $\{x_0\}$, then Condition \ref{A_sol unique} holds.
\end{theorem}
Theorem \ref{theorem: sol unique}  can be interpreted as follows: the original MC satisfies Condition \ref{A_sol unique} if it has a steady-state mass at $x_0\in\Omega$ such that the existence of stationary distribution as well as the existence of mass at $x_0$ is independent of the transitions from state $x_0$. 
In the rest of this subsection, we prove Theorem \ref{theorem: sol unique} and then give examples of independent mass singletons.

\subsubsection{Proof of Theorem \ref{theorem: sol unique}.}
We use the following two lemmas in our proof for Theorem \ref{theorem: sol unique}. 
\begin{lemma}[Balance Equation Reformulation II]\label{lemma:IE reform II}
%Consider an original MC in Definition \ref{A:Stationary Dist Eq} with stationary distribution $p$ and kernel $\bar{\kappa}$ satisfying Condition \ref{A_UBV}. 
Given the settings and assumptions in Theorem \ref{theorem: sol unique}, if constant $q :=\int_{\{x_0\}} dp(u)>0$, then $f=\frac{p}{q}$ is a solution to% the following integral equation:
\begin{align}
     & f(x) = \xi_0(x) + \int_{\Omega} \kappa_0(x,u) df(u), \quad \forall x\in \Omega,  \tag{BE-2$'$} \label{eqn:MC_IE_original prime}
\end{align}
where the inhomogeneous term $\xi_0\in\mathbf{\bar{X}}$ and the kernel $\kappa_0\in{T}(\Omega, \mathcal{B})$ are specified as
\begin{align*}
    & \xi_0(x):=\bar{\kappa}(x,x_0),\, x\in\Omega; \qquad \kappa_0(x,u) := \mathbf{1}\{u \neq x_0 \}\cdot\bar{\kappa}(x,u), \, x,u\in\Omega.
\end{align*}
Equivalently \eqref{eqn:MC_IE_original prime} can be stated as %an operator equation on space $\bar{\mathbf{X}}$:
\begin{align}
    & f = \xi_0+\mathcal{K}_0f, \tag{BE-3$'$} \label{eqn:MC_OE_original prime}
\end{align}
where operator $\mathcal{K}_0$ is defined by Definition \ref{Def:K} with kernel $\kappa_0$. 
Moreover, the solution space of \eqref{eqn:MC_BE_original} in $\bar{\mathbf{X}}$ is \rm{span}$(S)$ if $S\neq \varnothing$, where $S$ is the solution set of \eqref{eqn:MC_IE_original prime} or equivalently \eqref{eqn:MC_OE_original prime} in $\bar{\mathbf{X}}$. 
\end{lemma}
\begin{proof}{Proof:}
As a stationary distribution for the original MC, $p$ satisfies
\begin{align*}
    p(x) = \int_\Omega \bar{\kappa}(x,u) dp(u) = \int_{\{x_0\}} \bar{\kappa}(x,u) dp(u) + \int_{\Omega \backslash \{x_0\}} \bar{\kappa}(x,u) dp(u)=q\xi_0(x)+ \int_{\Omega} {\kappa}_0(x,u) dp(u),   
\end{align*}
for all $x\in\Omega.$ Multiplying both sides by $\frac{1}{q}$ yields that $f=\frac{p}{q}$ is a solution to \eqref{eqn:MC_IE_original prime}. 

With Definition \ref{Def:K}, we directly have the operator form of \eqref{eqn:MC_IE_original prime}, i.e., \eqref{eqn:MC_OE_original prime}. %By the definition, $\kappa_0$ satisfies Condition \ref{A_UBV} due to $\bar{\kappa}$ satisfying Condition \ref{A_UBV}. Then $\mathcal{K}_0$ is a well-defined continuous linear operator (Theorem 3.1 in \cite{li2021numerical}). Thus, equation \eqref{eqn:MC_OE_original prime} is well-defined. 

Now we only need to prove that the solution space of \eqref{eqn:MC_BE_original} in $\bar{\mathbf{X}}$ is \rm{span}$(S)$. Let $\bar{S}$ be the solution space of \eqref{eqn:MC_BE_original}. Consider any $f_1\in S$. We have that $f_1(x)=\xi_0(x)+\int_\Omega \kappa_0(x,u)df_1(u),\forall x\in\Omega$. Letting $x=1$ yields $f_1(1)=1 + f_1(1) - [f_1(x_0)-f_1(x_0-)] $. Thus, $f_1(x_0)-f_1(x_0-)=1$. Then $f_1(x)=\xi_0(x)+\int_\Omega \kappa_0(x,u)df_1(u)=\xi_0(x)[f_1(x_0)-f_1(x_0-)]+\int_\Omega \kappa_0(x,u)df_1(u)= \int_\Omega \mathbf{1}\{u=x_0\} \bar{\kappa}(x,u) df_1(u) +  \int_\Omega \mathbf{1}\{u \neq x_0\} \bar{\kappa}(x,u) df_1(u)=\int_\Omega \bar{\kappa}(x,u)df_1(u),\forall x\in\Omega$. Thus, $f_1\in\bar{S}$. Since $\bar{S}$ is a linear space, span$(S)\subseteq \bar{S}$.  

Next we use contradictions to show $\bar{S}\subseteq$span$(S)$. If $\bar{S}$ is not a subset to span$(S)$, there exists $f_2\in\bar{\mathbf{X}}$ such that $(i)$ $f_2\in \bar{S}$, $(ii)$ $f_2\notin$span$(S)$, $f_2\neq f_1$, $f_2\neq 0$, and $(iii)$ $f_2(x_0)-f_2(x_0-)\neq 0$. Here $(i)$ and $(ii)$ are straightforward. With regards to  $(iii)$, if $f_2(x_0)-f_2(x_0-)= 0$, we can construct a new $f'_2$ as $f_1+f_2$, which will satisfy $(i)$-$(iii)$. Then we have $f_2(x)=\int_\Omega \bar{\kappa}(x,u)df_2(u)= \int_\Omega \mathbf{1}\{u=x_0\} \bar{\kappa}(x,u) df_2(u) +  \int_\Omega \mathbf{1}\{u \neq x_0\} \bar{\kappa}(x,u) df_2(u)=\xi_0(x)[f_2(x_0)-f_2(x_0-)]+\int_\Omega \kappa_0(x,u)df_2(u),\forall x\in\Omega$. Multiplying both sides by $\frac{1}{f_2(x_0)-f_2(x_0-)}$ yields that $f_3(x)=\xi_0(x)+\int_\Omega \kappa_0(x,u) df_3(u),\forall x\in\Omega$, 
where $f_3:= \frac{f_2}{f_2(x_0)-f_2(x_0-)}$. Therefore, $f_3$ is also a solution to \eqref{eqn:MC_IE_original prime} and  $f_3\in S$. Because $f_2=f_3\cdot [f_2(x_0)-f_2(x_0-)]$, we have that $f_2\in$span$(S)$, which contradicts $f_2\notin$span$(S)$. Thus,  $\bar{S}\subseteq$span$(S)$. 
\Halmos
\end{proof} 
\begin{lemma}[Solution Uniqueness by Compactness]
\label{lemma:sol unique compactness}
%For an original MC in Definition \ref{A:Stationary Dist Eq}, if kernel $\bar{\kappa}$ satisfies Condition \ref{A_UBV}-\ref{A_MV}, 
Given the settings and assumptions in Theorem \ref{theorem: sol unique},  \eqref{eqn:MC_OE_original prime} has a unique solution in $\bar{\mathbf{X}}$ iff $\forall h\in\mathbf{D}$, $f=h+\mathcal{K}_0f$ has a solution $f\in\bar{\mathbf{X}}$.
\end{lemma}
\begin{proof}{Proof:}
%We first show that operator $\mathcal{K}'$ is compact. 
By the definition of $\kappa_0$, $\kappa_0$ satisfies Condition \ref{A_UBV}-\ref{A_MV} as a result of $\bar{\kappa}$ satisfying Condition \ref{A_UBV}-\ref{A_MV}. 
Then operator $\mathcal{K}_0$ is a well-defined compact linear operator (Theorem 3.1 and 3.4 in \cite{li2021numerical}). %Thus  \eqref{eqn:MC_OE_original} is a well-defined equation on $\bar{\mathbf{X}}$. 
Then according to the Riesz-Schauder's theorem (Lemma 2.4 in \cite{li2021numerical}), $\rm{Ker}(\mathcal{I}-\mathcal{K}_0)=\{0\}$ iff $\rm{Im}(\mathcal{I}-\mathcal{K}_0)=\bar{\mathbf{X}}$. We use this fact to prove Lemma \ref{lemma:sol unique compactness}.

($\Leftarrow$) Because for all $h\in\mathbf{D}$, there exists an $f\in\bar{\mathbf{X}}$ such that $f=h+\mathcal{K}_0f$, we have that 
${\mathbf{D}}\subseteq \rm{Im}(\mathcal{I}-\mathcal{K}_0)$. 
Because $\rm{Im}(\mathcal{I}-\mathcal{K}_0)$ is a liner space, ${\mathbf{X}}\subseteq \rm{Im}(\mathcal{I}-\mathcal{K}_0)$. 
Because $\mathcal{K}_0$ is a compact linear operator, $\rm{Im}(\mathcal{I}-\mathcal{K}_0)$ is closed. 
Thus, $\rm{Im}(\mathcal{I}-\mathcal{K}_0)=\bar{\mathbf{X}}$ and $\rm{Ker}(\mathcal{I}-\mathcal{K}_0)=\{0\}$. 
Thus, equation $f=\xi_0+\mathcal{K}_0f$, i.e., \eqref{eqn:MC_OE_original prime}, has a unique solution in $\bar{\mathbf{X}}$. 

($\Rightarrow$) Because equation $f=\xi_0+\mathcal{K}_0f$ has a unique solution, we have that $\rm{Ker}(\mathcal{I}-\mathcal{K}_0)=\{0\}$. Therefore, $\rm{Im}(\mathcal{I}-\mathcal{K}_0)=\bar{\mathbf{X}}$. Then for all $h\in\mathbf{D}$, there exists an $f\in\bar{\mathbf{X}}$ such that $f=h+\mathcal{K}_0f$.
\Halmos
\end{proof}

\begin{proof}{Proof of Theorem \ref{theorem: sol unique}:}
We first apply Lemma \ref{lemma:IE reform II}. 
According to Definition \ref{def:independent mass singleton} for independent mass singleton $\{x_0\}$, when $h= \bar \kappa(x,x_0)$, the transition kernel remains $\bar{\kappa}$ and has a stationary distribution $p$ with a mass at $x_0$. 
Then we obtain $q=\int_{\{x_0\}} dp(u)>0$. 
With all preconditions satisfied in Lemma \ref{lemma:IE reform II}, we obtain equation %$f=\xi_0+\mathcal{K}_0f$, i.e.,
\eqref{eqn:MC_OE_original prime} with a solution $f=\frac{p}{q}$. 

Next we consider an arbitrary $h\in\mathbf{D}$ replacing $\bar \kappa(\cdot,x_0)$. 
For the new MC with the modified kernel $\bar{\kappa}'$ defined in \eqref{eqn:kernel modified}, we have a stationary distribution $p'$ and a constant $q':=\int_{\{x_0\}}dp'(u)>0$  due to Definition \ref{def:independent mass singleton}. % for independent mass singleton $\{x_0\}$. 
By applying Lemma \ref{lemma:IE reform II} to the new MC with the modified kernel $\bar{\kappa}'$, we obtain equation $f=h+\mathcal{K}_0f$ with a solution $\frac{p'}{q'}$. Thus, $f=h+\mathcal{K}_0f$ has a solution $f\in\bar{\mathbf{X}}$ for all  $h\in\mathbf{D}$. 
This result exactly satisfies the requirement of Lemma \ref{lemma:sol unique compactness}. 
Then \eqref{eqn:MC_OE_original prime} has a unique solution. We have just shown $f=\frac{p}{q}$ is a solution and thus the unique solution. 
By lemma \ref{lemma:IE reform II}, the solution space of balance equation \eqref{eqn:MC_BE_original} is $\{\alpha p | \alpha \in\mathbb{R}\}$, and the original MC has a unique stationary distribution.
\begin{comment}
Finally, we only need to show that the original MC's stationary distribution is unique to complete the proof of Condition \ref{A_sol unique}. Suppose $p''\in\mathbf{D}$ is a stationary distribution. Then constant $q'':=\int_{\{x_0\}}dp''(u)>0$   due to Definition \ref{def:independent mass singleton}. Applying Lemma \ref{lemma:IE reform II} to the original MC yields that $f=\frac{p''}{q''}$ is a solution to \eqref{eqn:MC_OE_original prime}. Because \eqref{eqn:MC_OE_original prime} has a unique solution, we have $\frac{p}{q}=\frac{p''}{q''}$. Letting $x=1$ yields $\frac{1}{p}=\frac{1}{p''}$. Thus, $p=p''$, i.e., the original MC's stationary distribution is unique.
\end{comment}
\Halmos
\end{proof}

\subsubsection{Examples of independent mass singletons.}
We give examples of independent mass singletons based on accessibility and return time.
\begin{lemma}[Independent Mass Singleton Induced by Accessibility]\label{lemma:example of mass singleton 1}
In Approximation \ref{A:Stationary Dist Eq}, $x_0\in\Omega$ forms an independent mass singleton if $\exists m\in\mathbb{N}$ and $\delta>0$ such that the probability of returning to $x_0$ within $m$ transitions from any $u\in\Omega\backslash\{x_0\}$ is lower bounded by $\delta$, i.e., $\mathbb{P}_u [\tau_{\{x_0\}}\leqslant m]\geqslant \delta$, $\forall u\in\Omega\backslash\{x_0\} $. 
Here $\tau_{A}$ represents the first return time to $A\in\mathcal{B}$ and $\mathbb{P}_u$ represents the probability associated with the original MC's kernel $\bar{\kappa}$ and initial state $u$.
\end{lemma}
\begin{proof}{Proof:}
By Definition \ref{def:independent mass singleton}, we need to show that for all $h_1\in\mathbf{D}$, if we replace the definition of $\bar{\kappa}(\cdot,x_0)$ with $h_1$, the modified MC has a stationary  distribution and a steady-state mass at $x_0$.

To prove the existence of a stationary distribution, we will show that the modified MC is $\psi$-irreducible and recurrent. Moreover, $\mathbb{E}_{x_0}[\tau_{\{x_0\}}]$, the mean first return time from  $x_0$ to $x_0$ is finite. Then a stationary  distribution exists by Theorem 10.0.1 in \cite{meyn2012markov}.

(\textit{$\psi$-irreducible.})\quad Consider the modified MC. 
With any initial state $u\in\Omega\backslash\{x_0\}$, we have that $\tau_{\{x_0\}}$, i.e., the first return time to $x_0$, is independent of $h_1$, i.e., the transition from $x_0$ to $\Omega$.  
Therefore, we still have $\mathbb{P}_u [\tau_{\{x_0\}}\leqslant m]\geqslant \delta$ for all $u\in\Omega\backslash\{x_0\}$ under the modified MC.
This implies $x_0$ is accessible for all $x\in\Omega\backslash\{x_0\}$, and thus accessible for all $x\in\Omega$. 
Define measure $\varphi$: $\varphi(\{x_0\})=1$, $\varphi(\Omega\backslash\{x_0\})=0$. 
Then the modified MC is $\varphi$-irreducible and thus $\psi$-irreducible.

(\textit{Recurrent.}) \quad Because $\mathbb{P}_u [\tau_{\{x_0\}}\leqslant m]\geqslant \delta$ for all $u\in\Omega\backslash\{x_0\}$,  the probability of returning to $x_0$ within $m+1$ transitions from any state in $\Omega$ is lower bounded by $\delta$, i.e., $\mathbb{P}_u [\tau_{\{x_0\}}\leqslant m+1]\geqslant \delta$ for all $u\in\Omega$. By the strong Markov property, the  modified MC with an initial state $x_0$ returns to $\{x_0\}$ with probability 1, which implies the MC's recurrence by Theorem 8.3.6 in \cite{meyn2012markov}. 

(\textit{Finite mean first return time.})\quad The above argument also implies that the modified MC's mean first return time from $x_0$ to $x_0$  is upper bounded by $\frac{m+1}{\delta}$. Since $\{x_0\}$ as a singleton is a petite set, the modified MC has a stationary  distribution due to the above three properties ($\psi$-irreducibility, recurrence and finite mean first return time) by Theorem 10.0.1 in \cite{meyn2012markov}.  

Lastly, we prove the existence of mass at $x_0$. Let $\mathbf{\Phi}(t)(t\in\mathbb{N})$ be the modified MC and $p'$ be a stationary  distribution. Because $\mathbb{P}_u [\tau_{\{x_0\}}\leqslant m+1]\geqslant \delta$ for all $u\in\Omega$,  
we have that  $\max_{i=1,2,...,m+1}\mathbb{P}_u [\mathbf{\Phi}(i)=x_0]\geqslant \frac{\delta}{m+1}$ for all $u\in\Omega$. Thus, $\Omega= \cup_{i=1,2,...,m+1} C_i$, where $C_i:=\{u\in \Omega\,|\, \mathbb{P}_u [\mathbf{\Phi}(i)=x_0]\geqslant \frac{\delta}{m+1}\}$. There must exist $C_j$ $(j=1,2,...,m+1)$ such that $\int_{C_j}dp'(u)>0$. Because after $j$ transitions, the probability from any $u\in C_j$ to $x_0$ is lower bounded by $\frac{\delta}{m+1}$, we have that the steady-state probability of $\mathbf{\Phi}(t)= x_0$ is greater than $\frac{\delta}{m+1}\cdot \int_{C_j}dp'(u)>0$.
\Halmos
\end{proof}

\subsection{G/G/1+G Queue Example.}
\label{sec:GG1G}
%\subsubsection{Model Settings and Problem Formulation}
Consider a waiting list with a single server who follows a first-come-first-serve rule.
Arriving customers must wait in a buffer if the server is busy. 
A customer waiting in the buffer abandons after a random amount of patience time if her service has not been commenced, but never leaves during service. 
Once a customer being served completes service, the just-freed server immediately turns to the next customer in buffer as long as the buffer is not empty.
Customers' inter-arrival times, patience times and service times form three independent and identically distributed sequences that are independent of each other.
We aim to compute the following steady-state performance measures: the (customers') probability of no wait, the probability of abandonment, and queue length. 

For every customer, we define her virtual waiting time (VWT) as the amount of time she needs to wait before her service is commenced. 
Let $t_n$, $s_n$, $y_n$ and $w_n$ respectively be the inter-arrival time, service time, patience time and VWT of customer $n(n\in\mathbb{N})$.
Let $A(x)$ be the distribution of customers' inter-arrival times, $B(x)$ be the distribution of service times, and $G(x)$ be the distribution of patience times. 
We make the following general assumptions:
\begin{assumption}\label{A_GG1G}
The distributions associated with the G/G/1+G queue satisfy:
\begin{itemize}
    \item [$(i)$] (bounded support) $\exists \bar{y},\bar{s}\in\mathbb{R}_+$ such that $G(\bar{y})=B(\bar{s})=1$ and $\bar{y}+\bar{s}\leqslant 1$.  
    \item [$(ii)$] (smoothness) $\exists g_1,a_1,a_2\in\mathbb{R}_+$ such that $G'(x)\leqslant g_1, A'(x)\leqslant a_1, |A''(x)|\leqslant a_2, \forall x\in\Omega.$  
    \item [$(iii)$] (ergodicity) 
    $\mathbb{P}[s_n-\tau_{n+1}<0]>0$, i.e., the minimal service time is strictly smaller than the maximal inter-arrival time.
\end{itemize}
\end{assumption}
The assumption of $\bar{y}+\bar{s}\leqslant 1$ is without loss of generality and can be achieved by proper scaling. 
Item $(ii)$ assumes smoothness of the inter-arrival time distribution in terms of the first- and second-order derivatives. 
Example distributions include exponential distributions, Erlang distributions, Gamma distributions (with a shape parameter no smaller than 2), distributions that can be expressed as a polynomial, etc. 
Similarly, item $(ii)$ also assumes the smoothness of the patience time distribution in terms of the first-order derivative. 
Item $(iii)$ is exactly the same with the ergodicity condition for G/G/1+G queue provided by \citep{baccelli1981queues}.

Customers' VWTs are described using an MC: 
for all $n\in\mathbb{N}$, we have
\begin{align}
    w_{n+1} = \begin{cases}
        [w_n+s_{n}-t_{n+1}]_+, \quad \text{if\,} y_n > w_n, \\
        [w_n-t_{n+1}]_+, \quad \text{otherwise}.
    \end{cases} \label{eqn:gg1g original MC}
\end{align}
%Here $[\cdot]$ as an operation is defined as $[x]_+=\max\{x,0\}$.
Due to $(i)$ of Assumption \ref{A_GG1G}, the MC is supported on $\Omega=[0,1]$. 
Let $p(x)\in\mathbf{D}$ be the stationary distribution of $w_n$. 
The steady-state performance measures we aim to compute can be written as 
\begin{align}\label{eqn:gg1g sspm}
    \mathbb{E}\phi_i(X)= \int_\Omega \phi_i(x) dp(x), \quad i\in\{1,2,3\},
\end{align}
where $\phi_1(x)=\mathbf{1}\{x=0\}$ computes the probability of no wait, $\phi_2(x)=G(x)$ computes the probability of abandonment, and $\phi_3(x)=\lambda h(x)$ computes the queue length with $h(x):=\int_{\mathbb{R}_+} \min\{y,x\}dG(y)$ and $\lambda:={1}/{\int_{\mathbb{R}_+}x dA(x)} $.
The balance equation can be written as: $\forall x\in\Omega$,
\begin{align}\label{eqn:gg1g sde}
    p(x) = \int_\Omega \bar{\kappa}(x,u) dp(u),\quad \text{where } \bar{\kappa}(x,u)=\bar{G}(u)D(x-u)+G(u){A}^*(u-x).
\end{align}
Here $D(x)$ is the distribution of $s_n-\tau_{n+1}$ and ${A}^*(x):=1-{A}(x-)$. 
Thus, the MC in \eqref{eqn:gg1g original MC} defined by the VWTs is in accordance with the definition of original MC in Approximation \ref{A:Stationary Dist Eq}.

We apply Approximation \ref{A:Stationary Dist Eq} and generate a sequence of finite-state MCs to approximate the original MC described in \eqref{eqn:gg1g original MC}. 
The $r$-th MC has $J^{(r)}=2^r+1$ states.
Its states $\{c^{(r)}_{i}\}_{i=1}^{J^{(r)}}$ and  transition matrix $Q^{(r)}=\{q^{(r)}_{{ij}}\}_{i,j=1,2,...,J^{(r)}}$ are respectively defined as
\begin{align}
    c^{(r)}_{i} := \frac{(i-1)}{2^r},\, i=1,2,...,J^{(r)};\quad q_{{ij}}^{(r)}:= \bar{\kappa}(c^{(r)}_{j},c^{(r)}_{i})-\bar{\kappa}(c^{(r)}_{{j-1}},c^{(r)}_{i}),\, i,j=1,2,...,J^{(r)}.
    \label{eqn:gg1g approximate MC}
\end{align}
Particularly, $c^{(r)}_{0}\equiv -\infty$. 
Let $p^{(r)}(x)$ be the discrete stationary distribution of the $r$-th approximate MC. 
We verify Condition \ref{A_UBV}-\ref{A_T2} so that Theorem \ref{Theorem: convergence naive} is applicable.  
\begin{theorem}\label{theorem:gg1g conditions} 
Consider a G/G/1+G queue with VWTs described as an MC in \eqref{eqn:gg1g original MC}. The finite approximation outlined in \eqref{eqn:gg1g approximate MC} satisfies Condition \ref{A_UBV}-\ref{A_T2} under Assumption \ref{A_GG1G}.  
\end{theorem}
\begin{proof}{Proof:}
(\textit{Verification of Condition \ref{A_UBV}}.)\quad By definition \eqref{eqn:gg1g sde}, the kernel $\bar{\kappa}(x,u)$ can be written as 
\begin{align}
    \bar\kappa(x,u)=\int_\Omega [\bar{G}(u)A^*(u-x+s) + G(u)A^*(u-x)] dB(s), \quad \forall x,u\in\Omega.\label{eqn:gg1g kernel form2}
\end{align}
Then $|\frac{\partial \bar\kappa(x,u)}{\partial u}|\leqslant \int_\Omega | \frac{\partial \bar{G}(u)}{\partial u} A^*(u-x+s)+ \frac{\partial A^*(u-x+s)}{\partial u} \bar{G}(u)+ \frac{\partial  G(u)}{\partial u} A^*(u-x) + \frac{\partial A^*(u-x)}{\partial u} G(u)| d B(s)\leqslant \int_\Omega \big[ g_1 |A^*(u-x+s)-A^*(u-x)| + a_1 \bar{G}(u)  + a_1 G(u) \big] dB(s)\leqslant g_1+a_1$. As discussed in Subsection \ref{sec:verify c1-3}, $\bar{\kappa}$ satisfies Condition \ref{A_UBV}.  

(\textit{Verification of Condition \ref{A_MV}}.)\quad By  \eqref{eqn:gg1g kernel form2}, $|\frac{\partial^2 \bar\kappa(x,u)}{\partial x \partial u}|\leqslant \int_\Omega | \frac{\partial \bar{G}(u)}{\partial u} \frac{\partial A^*(u-x+s) }{\partial x}  + \frac{\partial^2 A^*(u-x+s)}{\partial x \partial u} \bar{G}(u)+ \frac{\partial G(u)}{\partial u} \frac{\partial A^*(u-x)}{\partial x} + \frac{\partial^2 A^*(u-x)}{\partial x \partial u} G(u)| dB(s)\leqslant \int_\Omega | g_1 a_1 + a_2 \bar{G}(u)+ a_1g_1  + a_2 {G}(u) | dB(s)\leqslant 2g_1 a_1 + a_2.$
In addition, $|\frac{d \bar{\kappa}(x,1)}{d x}|=|\frac{d A^*(1-x) }{d x}|\leqslant a_1$.
Then as discussed in Subsection \ref{sec:verify c1-3}, $\bar{\kappa}$ satisfies Condition \ref{A_MV}.

(\textit{Verification of Condition \ref{A_sol unique}.})\quad By Theorem \ref{theorem: sol unique}, we only need to show $\{0\}$ is an independent mass singleton.
Due to $(iii)$ of Assumption \ref{A_GG1G}, there exist $l,\delta>0$ such that $\mathbb{P}[s_n-\tau_{n+1}\leqslant -l]\geqslant \delta$. 
Then by \eqref{eqn:gg1g original MC}, we have $\mathbb{P}[w_{n+1} \leqslant \max\{w_n-l,0\}]\geqslant \mathbb{P}[s_n-\tau_{n+1}\leqslant -l]\geqslant \delta$.
Thus, $\mathbb{P}_u [\tau_{\{0\}}\leqslant \lceil \frac{1}{l} \rceil ]\geqslant \delta^{\lceil \frac{1}{l} \rceil}$ for all $u\in [0,1]\backslash\{0\}$. 
By Lemma \ref{lemma:example of mass singleton 1}, $\{0\}$ is an independent mass singleton. 

(\textit{Verification of Condition \ref{A_T2}.})\quad
As discussed in Subsection \ref{sec:verify c1-3}, the approximate MCs' states are in accordance with the finite splits in Condition \ref{A_MV}. 
Thus, $(i)$ of Condition \ref{A_T2} is satisfied. 
Item $(ii)$ and $(iii)$ of Condition \ref{A_T2} are directly satisfied by definition \eqref{eqn:gg1g approximate MC}. 
Item $(iv)$ of Condition \ref{A_T2} is guaranteed by the smoothness of $\bar{\kappa}$ with $\frac{\partial \bar \kappa(x,u)}{\partial u}$ and $\frac{\partial \bar \kappa(x,u)}{\partial x}$ both bounded for all $(x,u)\in \Omega^2$. Here the bound of $\frac{\partial \bar \kappa(x,u)}{\partial u}$ has been provided when verifying Condition \ref{A_MV}. 
With regards to the bound of $\frac{\partial \bar \kappa(x,u)}{\partial x}$, we have $|\frac{\partial \bar\kappa(x,u)}{\partial x}|\leqslant \int_\Omega | \frac{\partial A^*(u-x+s)}{\partial x} \bar{G}(u) + \frac{\partial A^*(u-x)}{\partial x} G(u)| d B(s)\leqslant \int_\Omega |a_1\bar{G}(u)+a_1{G}(u)| dB(s)\leqslant a_1.$
\Halmos
\end{proof}

%Additional examples that satisfy the independent mass singleton assumption with queuing settings can be found in our companion paper \cite{li2021modeling}.

%\subsection{Verify Condition 4 of proper truncation}

\section{Numerical Experiments.}\label{sec:Numerical}

\subsection{Settings.} 
We consider steady-state performance evaluation of an overloaded G/G/1+G queue. The steady-state performance measures we compute are given in \eqref{eqn:gg1g sspm}. The inter-arrival times are assumed to follow a scaled Erlang distribution: $A(x)=A_0(\lambda x)$, $A_0(x)\sim\text{Er}(2,2)$.
The service times follow a scaled Beta distribution: $B(x)=B_0(\frac{\mu x}{2}), B_0(x)\sim \text{Beta}(2,2)$.
The patience times also follow a scaled Beta distribution: $G(x)=G_0(2 x), G_0(x)\sim \text{Beta}(3,4)$.
Here $\lambda$ and $\mu$ are respectively the arrival intensity and service rate.  
We respectively compute for the models $(a)$ $(\lambda,\mu)=(4.1,4)$ and $(b)$ $(\lambda,\mu)=(5,4)$.

We compare our finite approximation approach against a standard MCMC implementation and a fluid approximation provided by \cite{whitt2006fluid}. 
Our finite approximation approach is as outlined in \eqref{eqn:gg1g approximate MC}. 
In MCMC, we use Gibbs sampler with a burn-in period of  $10^6$ steps, and then collect one sample for every $10^3$ steps. 
In fluid approximation, the stationary distribution of VWT is considered to be a one-point distribution concentrated at $w_{\text{Fluid}}=G^{-1}(1-\frac{\mu}{\lambda})$. 
We refer to the approximate MC's stationary distribution in finite approximation, the emprical distribution in MCMC, and the one-point distribution in fluid approximation as \textit{proxy stationary distributions} of the VWT because for all these three approaches, we are essentially using the proxy stationary distributions to compute steady-state performance measures. 
Moreover, the accuracy of the proxy stationary distribution decides that of steady-state performance evaluation. 
We compute and compare:
\begin{itemize}
    \item 
    Balance equation residuals. Because the true stationary distribution is unavailable, we are not able to directly compute errors of proxy stationary distributions. 
    Instead, we use the $L^\infty$ residues in balance equations \eqref{eqn:gg1g sde} to measure the proxy stationary distributions' accuracy. 
    We compute residues by discretizing support $[0,1]$ and enumerating $\{ x\in[0,1] | 10^7x\in\mathbb{N}\}$.  
    \item Computation time.
    \item Approximate/estimated steady-state performance measures. 
    Note that finite and fluid approximations yield deterministic approximators, while MCMC yields probabilistic estimators. 
    \item Error bounds/confidence intervals. Note that finite approximation yields strict deterministic bounds for steady-state performance measures, while MCMC yields confidence intervals. All confidence intervals are constructed via central limit theorem at the level of $99.9\%$. The fluid approximation approach does not have available computable error bounds.
\end{itemize}
%The computation time reflects efficiency of different approaches. As regards performance evaluation, we expect a superior approach to yield more accurate approximated/estimated steady-state performance measures, smaller error bounds/confidence intervals, and faster convergence rates.         

For both model $(a)$ and $(b)$, we run a series of experiments using finite approximation and MCMC. 
In the $r$-th experiment, finite approximation uses an approximate MC with $N=2^r+1$ states, while MCMC collects $N=2^r+1$ samples. 
In this way, the proxy stationary distributions obtained from finite approximation and MCMC are both discrete distributions with $N$ jumps, and $N$ can be regarded as a measure of distribution granularity. 
Fluid approximation is computed once. 

\subsection{Results.}  
\begin{figure}[!htb]
    \caption{A comparison of finite approximation vs. MCMC  with respect to $L^\infty$ residues in balance equations \eqref{eqn:gg1g sde}. ($N$ = distribution granularity, i.e., the number of jumps in the proxy stationary distributions.)}\label{fig:residue plot} 
    \centering
    \includegraphics[width=0.45\textwidth,height=100px]{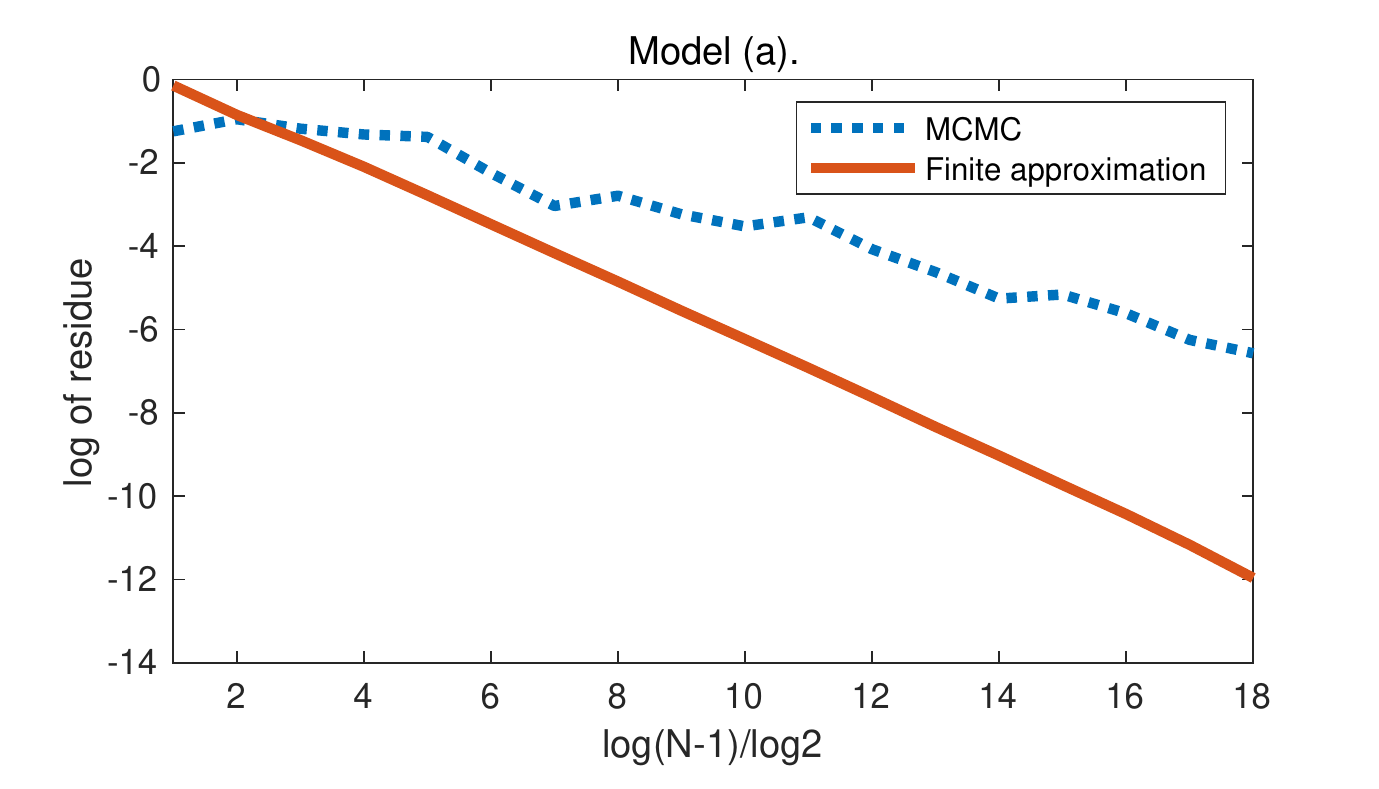} \qquad 
    \includegraphics[width=0.45\textwidth,height=100px]{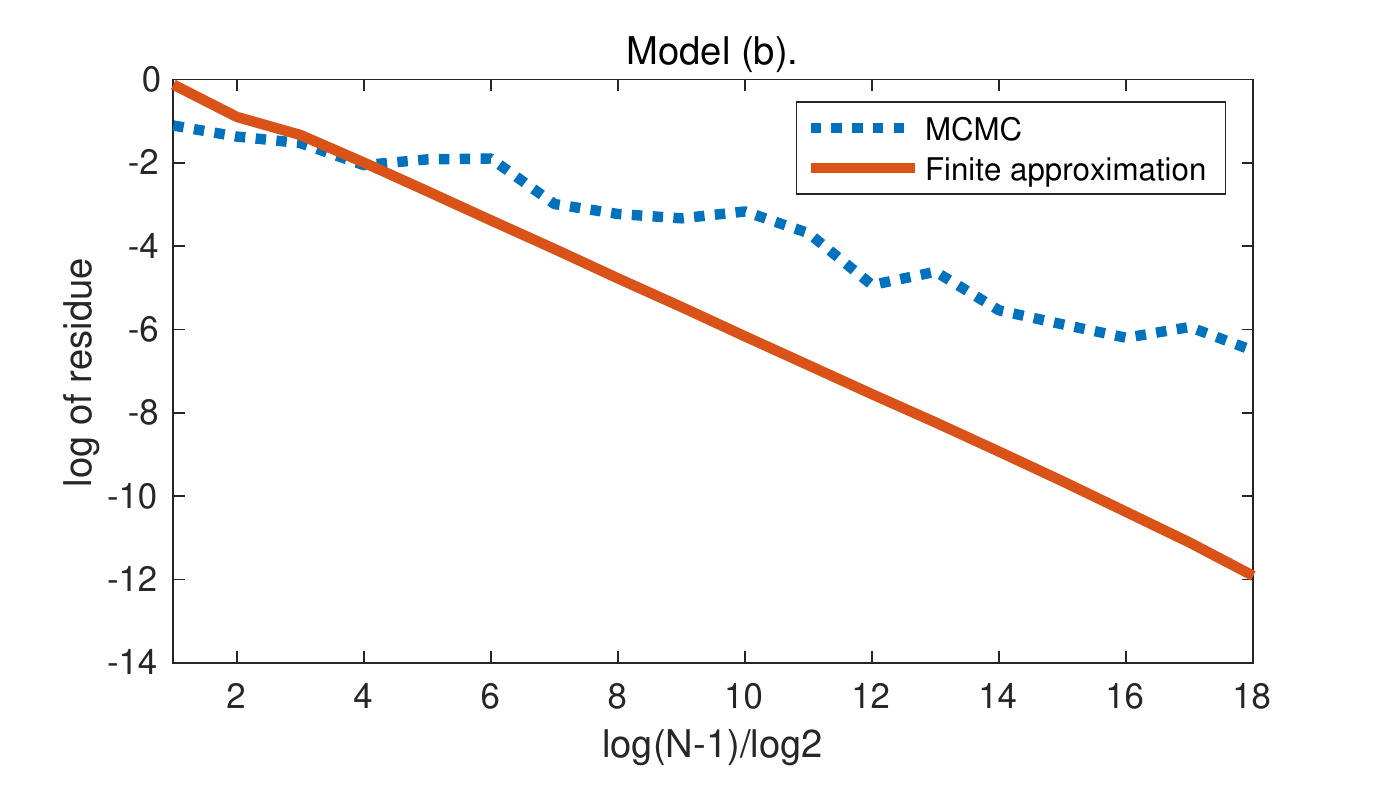}
\end{figure}

\,\\ 
\noindent \textbf{Finite Approximation vs. MCMC.} \quad 
We compare proxy stationary distributions' $L^\infty$ residues in balance equations \eqref{eqn:gg1g sde} (Figure \ref{fig:residue plot}): 
finite approximation has smaller residues and converges faster than MCMC. 
Although the residues are not necessarily linear in solution errors for balance equations \eqref{eqn:gg1g sde}, we find the finite approximation residues converging as $O(\frac{1}{N})$ approximately, which is in accordance with the convergence rate provided in Subsection \ref{sec:near optimality}. 
With regards to MCMC, we find the residues converging as $O(\frac{1}{\sqrt{N}})$ approximately, which is in accordance with the convergence rate provided by the DKW inequality \citep{kosorok2006introduction}. 
We also find that if we change the measure into $L^1$ residues, finite approximation still outperforms MCMC with similar advantages (results available from the authors).

Table \ref{table: residue stats} presents balance equation residues and corresponding computation times of different approaches. 
In comparison with MCMC, finite approximation needs significantly less time to achieve the same level of accuracy. 
For example, MCMC needs as many as $2^{18}+1\approx 260,000$ samples and $42$-minute computation time to reduce residues to the level of $1\times 10^{-3}$ in both Model $(a)$ and $(b)$. 
Finite approximation only needs to construct an MC of $2^9+1=513$  states to achieve the same level of residues, and the computation time is about $0.04$ second. 
This shows the dramatic improvement our approach provides over the use of MCMC. 
As $N$ increases, a sample of MCMC takes $O(N)$ time as expected. 
For finite approximation, our computation based on iterations takes $O(zN^2)$ time, where $z$ is the number of iterations. 
 This computation cost can be further reduced if the approximate finite-state MC has a factorization structure \cite{barreto2011computing}. 

\begin{table}[!htb]
\centering
\caption{A comparison of finite approximation vs. MCMC  with respect to $L^\infty$ residues in balance equations \eqref{eqn:gg1g sde}, and computation time. ($N$ = distribution granularity, i.e., the number of jumps in the proxy stationary distributions; CT = computation time/second.)}
\label{table: residue stats}
\tiny 
\begin{tabular}{ccccccccccc}
\multicolumn{5}{c}{Model $(a)$}                                                                                                                &                       & \multicolumn{5}{c}{Model $(b)$}                                                                                                               \\ \cline{1-5} \cline{7-11} 
\multicolumn{1}{|c|}{\multirow{2}{*}{$\log_2(N-1)$}} & \multicolumn{2}{c|}{Finite approximation}    & \multicolumn{2}{c|}{MCMC}                    & \multicolumn{1}{c|}{} & \multicolumn{1}{c|}{\multirow{2}{*}{$\log_2(N-1)$}} & \multicolumn{2}{c|}{Finite approximation}    & \multicolumn{2}{c|}{MCMC}                    \\
\multicolumn{1}{|c|}{}                       & $L^\infty$ residue & \multicolumn{1}{c|}{CT}       & $L^\infty$ residue & \multicolumn{1}{c|}{CT}       & \multicolumn{1}{c|}{} & \multicolumn{1}{c|}{}                       & $L^\infty$ residue & \multicolumn{1}{c|}{CT}       & $L^\infty$ residue & \multicolumn{1}{c|}{CT}       \\ \cline{1-5} \cline{7-11} 
\multicolumn{1}{|c|}{3}                  & 2.33E-01  & \multicolumn{1}{c|}{1.56E-03} & 3.07E-01 & \multicolumn{1}{c|}{9.83E+00} & \multicolumn{1}{c|}{} & \multicolumn{1}{c|}{3}                  & 2.64E-01  & \multicolumn{1}{c|}{1.69E-03} & 2.17E-01 & \multicolumn{1}{c|}{9.53E+00} \\
\multicolumn{1}{|c|}{6}                  & 3.10E-02  & \multicolumn{1}{c|}{5.70E-03} & 1.06E-01 & \multicolumn{1}{c|}{1.04E+01} & \multicolumn{1}{c|}{} & \multicolumn{1}{c|}{6}                  & 3.40E-02  & \multicolumn{1}{c|}{6.04E-03} & 1.50E-01 & \multicolumn{1}{c|}{1.01E+01} \\
\multicolumn{1}{|c|}{9}                  & 3.88E-03  & \multicolumn{1}{c|}{4.03E-02} & 3.93E-02 & \multicolumn{1}{c|}{1.47E+01} & \multicolumn{1}{c|}{} & \multicolumn{1}{c|}{9}                  & 4.25E-03  & \multicolumn{1}{c|}{4.10E-02} & 3.57E-02 & \multicolumn{1}{c|}{1.42E+01} \\
\multicolumn{1}{|c|}{12}                 & 4.88E-04  & \multicolumn{1}{c|}{7.66E-01} & 1.70E-02 & \multicolumn{1}{c|}{4.95E+01} & \multicolumn{1}{c|}{} & \multicolumn{1}{c|}{12}                 & 5.24E-04  & \multicolumn{1}{c|}{8.11E-01} & 7.20E-03 & \multicolumn{1}{c|}{4.85E+01} \\
\multicolumn{1}{|c|}{15}                 & 5.93E-05  & \multicolumn{1}{c|}{3.25E+01} & 5.74E-03 & \multicolumn{1}{c|}{3.31E+02} & \multicolumn{1}{c|}{} & \multicolumn{1}{c|}{15}                 & 6.48E-05  & \multicolumn{1}{c|}{3.69E+01} & 2.80E-03 & \multicolumn{1}{c|}{3.27E+02} \\
\multicolumn{1}{|c|}{18}                 & 6.37E-06  & \multicolumn{1}{c|}{2.31E+03} & 1.40E-03 & \multicolumn{1}{c|}{2.53E+03} & \multicolumn{1}{c|}{} & \multicolumn{1}{c|}{18}                 & 6.70E-06  & \multicolumn{1}{c|}{2.49E+03} & 1.49E-03 & \multicolumn{1}{c|}{2.53E+03} \\ \cline{1-5}
\cline{7-11}
\end{tabular}
\end{table}

\begin{figure}[!htb]
    \caption{A comparison of finite approximation vs. MCMC with respect to  steady-state performance evaluation, results from Model $(a)$. ($N$ = distribution granularity, i.e., the number of jumps included in the proxy stationary distributions. Note that finite approximation provides approximate steady-state performance measures with bounds, while MCMC  provides estimated steady-state performance measures with  99.9\% confidence intervals.)}\label{fig:measure plot} 
    \centering
    \tiny
    Evaluate $\mathbb{E}\phi_1(X)= \mathbb{E}\mathbf{1}\{X=0\}$ in Model $(a)$.
    \includegraphics[width=1\textwidth, height = 85px]{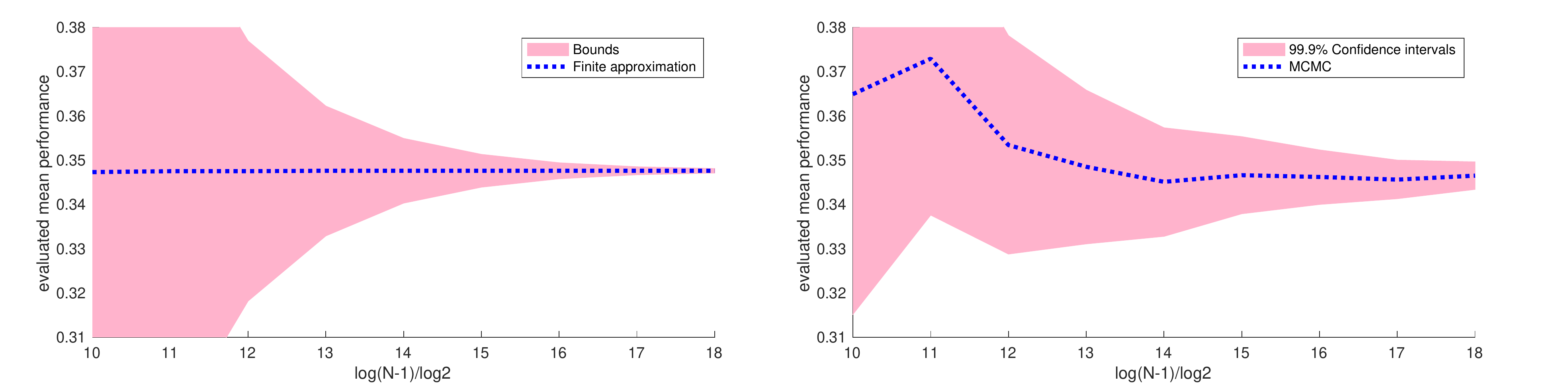} \\ 
    Evaluate $\mathbb{E}\phi_2(X)= \mathbb{E}G(X)$ in Model $(a)$.
    \includegraphics[width=1\textwidth, height = 85px]{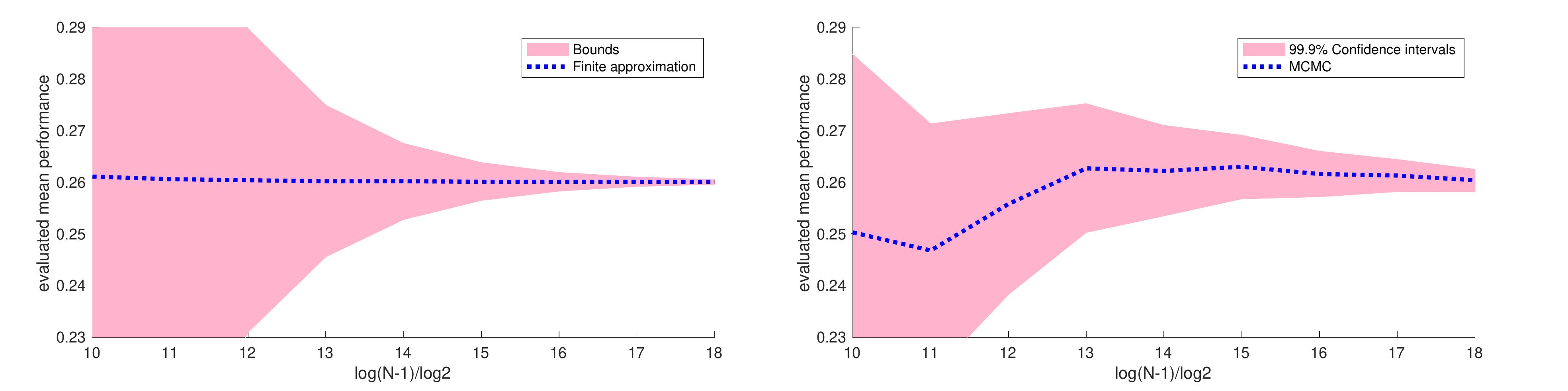} \\
    Evaluate $\mathbb{E}\phi_3(X)= \mathbb{E}\lambda h(X)$ in Model $(a)$.
    \includegraphics[width=1\textwidth, height = 85px]{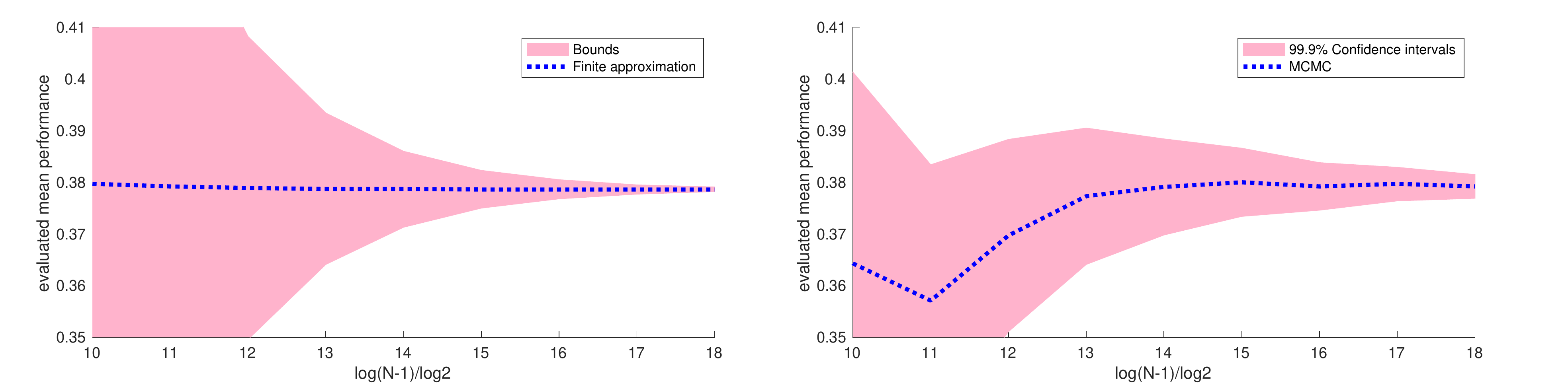} 
 \end{figure}

\begin{figure}[!htb]
    \caption{A comparison of finite approximation vs. MCMC with respect to  steady-state performance evaluation, results from Model $(b)$. ($N$ = distribution granularity, i.e., the number of jumps included in the proxy stationary distribution. Note that finite approximation provides approximate steady-state performance measures with bounds, while MCMC  provides estimated steady-state performance measures with  99.9\% confidence intervals.)}\label{fig:measure plot 2} 
    \centering
    \tiny
    Evaluate $\mathbb{E}\phi_1(X)= \mathbb{E}\mathbf{1}\{X=0\}$ in Model $(b)$.
    \includegraphics[width=1\textwidth, height = 85px]{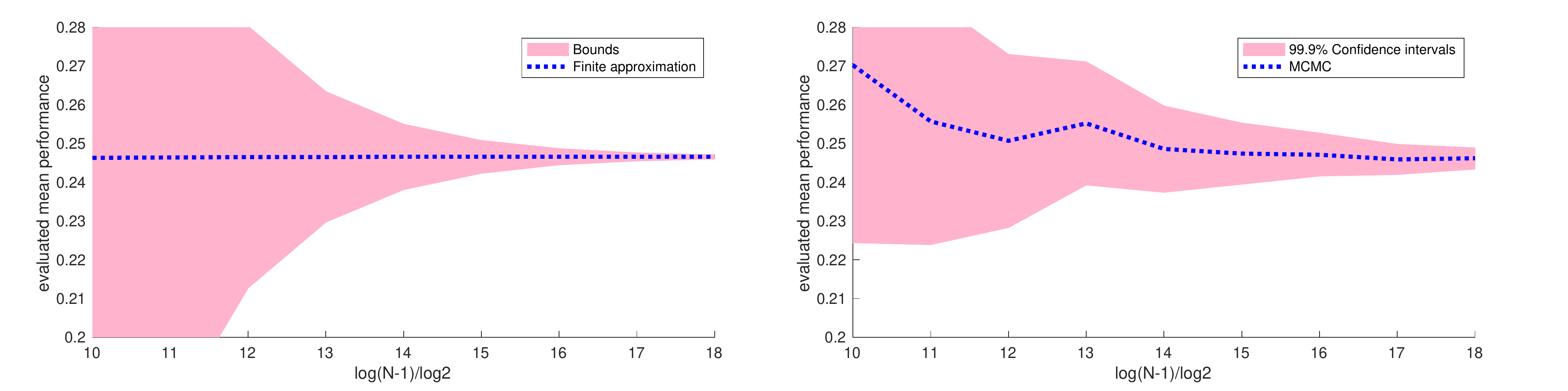} \\ 
    Evaluate $\mathbb{E}\phi_2(X)= \mathbb{E}G(X)$ in Model $(b)$.
    \includegraphics[width=1\textwidth, height = 85px]{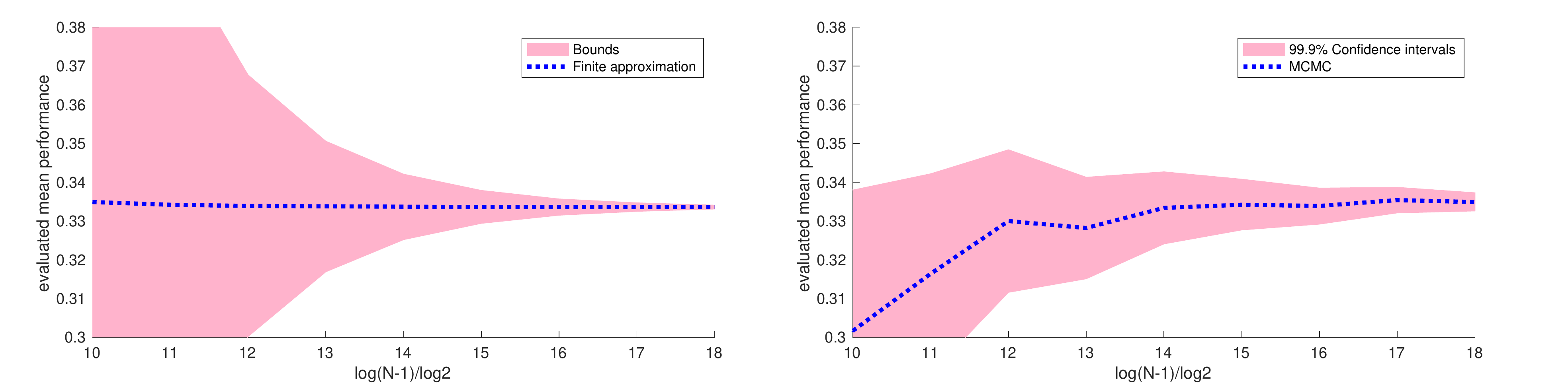} \\
    Evaluate $\mathbb{E}\phi_3(X)= \mathbb{E}\lambda h(X)$ in Model $(b)$.
    \includegraphics[width=1\textwidth, height = 85px]{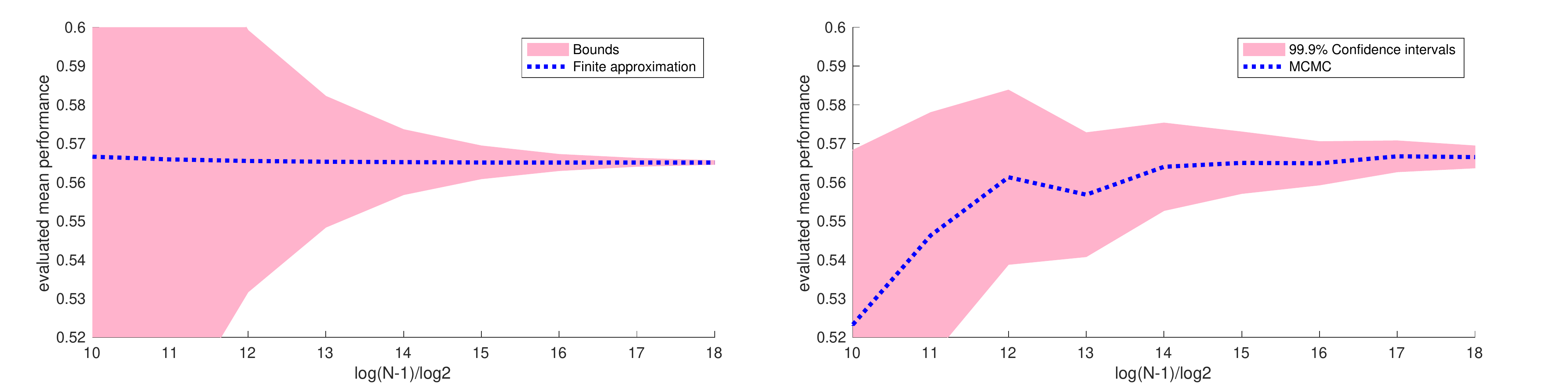} 
 \end{figure}

Lastly, we compare the accuracy and efficiency in performance evaluation (Figure \ref{fig:measure plot}-\ref{fig:measure plot 2}).
%Note that in Model $(b)$, the BRW model is ``unbalanced'' due to $\mathbb{E}W_t>0$: its stationary distribution is concentrated around $1$ and has a tiny variance.  Thus, Model $(b)$ favors MCMC in performance evaluation. Nevertheless, even in Model $(b)$, finite approximation outperforms MCMC as observed below.

\begin{itemize}
    \item Accuracy: finite approximators have significantly smaller errors than MCMC estimators. 
    For example, when we compute $\mathbb{E}\phi_1(X)$, or the probability of no wait, MCMC still has an error at the level $0.1\%$ with as many as 260,000 samples and $42$-minute computation time in both Models. 
    The finite approximator achieves the same level of accuracy by constructing an MC of $2^{8}+1=257$ states and the computation time is around $0.02$ second. 
    Here relative errors are computed by using the average results of finite approximation ($N=18$) and MCMC $(N=21)$ as reference values.
    %Moreover, by increasing the MC state count to $2^{18}+1\approx 260,000$ and using around $40$ minutes' computation time, finite approximation is able to reduce the error to $0.1\%$. %Results for other performance measure functions are similar. 
    \item Efficiency: The performance measure bounds provided by finite approximation are converging faster than the confidence intervals provided by MCMC. 
    Figures \ref{fig:measure plot}-\ref{fig:measure plot 2} indicate that the finite approximation bound length is converging as $O(\frac{1}{N})$, which is in accordance with Theorem \ref{theorem:performance measure error} and the convergence rate provided in Subsection \ref{sec:near optimality}. 
    By comparison, the MCMC confidence interval length is converging as $O(\frac{1}{\sqrt{N}})$, which is in accordance with the central limit theorem. 
    When $N=2^{18}+1\approx 260,000$, the finite approximation bound length has reduced to $10^{-3}$, while the  MCMC confidence interval length is still at the level of $10^{-2}$.
\end{itemize}

\noindent \textbf{Finite Approximation vs. Fluid Approximation.} \quad
Because the large market assumption does not hold in our experiments, the fluid approximation has huge errors (Table \ref{table:FA vs FA}). 
%: in Model $(a)$, the $L^\infty$ residue for balance equations is 0.544, the $L^1$ residue is 0.106, $\mathbb{E}\phi_1(X)=0$, $\mathbb{E}\phi_2(X)=0.024$ and $\mathbb{E}\phi_3(X)=0.239$. 
%In Model $(b)$, the $L^\infty$ residue for balance equations is 0.507, the $L^1$ residue is 0.129, $\mathbb{E}\phi_1(X)=0$, $\mathbb{E}\phi_2(X)=0.2$ and $\mathbb{E}\phi_3(X)=0.64$.
The results show the dramatic advantage of our approach in analyzing limited-scale stochastic systems.

\begin{table}[!htb]
\caption{A comparison of finite approximation vs. fluid approximation  with respect to balance equation residues and steady-state performance measures. (Prob.= Probability. $N=2^{18}+1$ is used for finite approximation. *: relative errors computed with MCMC results $(N=21)$ as reference values.)}
\centering
\label{table:FA vs FA}
\tiny
\begin{tabular}{ccccccc|}
\multicolumn{3}{c}{Model $(a)$}                                                                                                                                                               &                       & \multicolumn{3}{c}{Model $(b)$}                                                                                                                                                              \\ \cline{1-3} \cline{5-7} 
\multicolumn{1}{|c|}{}  & \begin{tabular}[c]{@{}c@{}}Finite \\ Approximation\end{tabular} & \begin{tabular}[c]{@{}c@{}}Fluid\\ Approximation\end{tabular}     & \multicolumn{1}{|c|}{} & \multicolumn{1}{c|}{}  & \begin{tabular}[c]{@{}c@{}}Finite \\ Approximation\end{tabular} & \begin{tabular}[c|]{@{}c@{}}Fluid\\ Approximation\end{tabular}    \\ \cline{1-3} \cline{5-7} 
\multicolumn{1}{|c|}{$L^\infty$ residues} & 6.37E-06                                                        & 5.44E-01                                                  & \multicolumn{1}{|c|}{} & \multicolumn{1}{c|}{$L^\infty$ residues} & 6.70E-06                                                        & 5.07E-01                                                  \\
\multicolumn{1}{|c|}{$L^1$ residues} & 8.68E-07                                                        & 1.06E-01                                                    & \multicolumn{1}{|c|}{} & \multicolumn{1}{c|}{$L^1$ residues} & 9.69E-07                                                        & 1.29E-01                                                   \\
\multicolumn{1}{|c|}{Prob. of no Wait$^*$} & 0.0005\%                                                        & 100\%                                                      & \multicolumn{1}{|c|}{} & \multicolumn{1}{c|}{Prob. of no Wait$^*$} & 0.0122\%                                                          & 100\%                                                     \\
\multicolumn{1}{|c|}{Prob. of Abandonment$^*$} & 0.0035\%                                                        & 91\%                                             & \multicolumn{1}{|c|}{} & \multicolumn{1}{c|}{Prob. of Abandonment$^*$} & 0.0004\%                                                         & 40\%                                               \\
\multicolumn{1}{|c|}{Queue Length$^*$} & 0.0027\%                                                         & 37\%                                           & \multicolumn{1}{|c|}{} & \multicolumn{1}{c|}{Queue Length$^*$} & 0.0001\%                                                        & 13\%                                                    \\ \cline{1-3} \cline{5-7} 
\end{tabular}
\end{table}

\section{Conclusion.}
\label{sec:ConcludingRemarks}
This paper develops a general and consistent finite approximation scheme with tractable solutions and error bounds for computing the stationary distribution and associated steady-state performance measures of an MC supported on $\mathbb{R}$. 
The approximate MCs are constructed based on variation properties of the original MC's transition kernel. 
Our approach is also applicable for countable-state MCs, and those mixing countable states with continuous states. We have also shown that our approach is near-optimal among all approximation methods using discrete distributions, given a Lipschitz continuity assumption on the original MC’s transition kernel. In terms of its rate of convergence, our approach is theoretically more efficient than using MCMC. 
Another possible reason for our approach's advantage over MCMC is that in MCMC's proxy stationary distribution (i.e., the empirical distribution), states are random and the jump sizes are uniform and not weighted. In contrast, finite approximation selects states and weighs jump sizes via constructing an appropriate transition matrix.

A series of numerical experiments show that the method developed in this paper outperforms benchmark results when compared with MCMC and fluid approximation in its efficiency and accuracy by several orders of magnitude. These results also show a 
major advantage of the approach developed in  this paper in computing performance measures of limited-scale stochastic systems.

%We close this paper by pointing out three plausible directions for future research. $(a)$ This paper studies approximation based on finite-state MCs. As a result, the approximate stationary distribution is discrete and non-smooth. An interesting direction is to reduce approximation errors by increasing the smoothness of approximate stationary distributions (e.g., via extra iterations). $(b)$ This paper considers MCs supported on $\mathbb{R}$. It is important to extend our analysis methods for multi-dimensional MCs.$(c)$ This paper only considers the performance evaluation problems. An interesting direction is to integrate our techniques with existing optimization frameworks (e.g., the response surface methodology), develop optimization algorithms, and apply to decision problems. 

\bigskip

% Appendix here
% Options are (1) APPENDIX (with or without general title) or 
%             (2) APPENDICES (if it has more than one unrelated sections)
% Outcomment the appropriate case if necessary
%
\begin{APPENDICES}
\section{List of Cited Theorems and Lemmas from the Literature.}
\label{append:cited theorems}
\begin{theorem}[Theorem 2.5 in \cite{li2021numerical}]
Let $\bar{\mathbf{X}}$ be the closure of  $\mathbf{X}$ under the norm $||\cdot||_\infty$. Then $(\bar{\mathbf{X}}, ||\cdot||_\infty)$ is a Banach space. 
\end{theorem}

\begin{theorem}[Theorem 3.1 in \cite{li2021numerical}]
For $\tau \in {T}(\Omega,\mathcal{B})$, operator $\mathcal{L}$ in Definition \ref{Def:K} is a well defined  continuous linear operator on $\bar{\mathbf{X}}$ if  Condition \ref{A_UBV} (uniformly bounded variation) holds. Conversely, if a continuous linear operator $\mathcal{L}$ on $\bar{\mathbf{X}}$ has the form \eqref{Def:K_1}, then Condition \ref{A_UBV} must hold. 
Particularly, 
\begin{align*}
    ||\mathcal{L}||_O\leqslant  \sup_{x\in\Omega} [2V_u(\tau(x,u))+ \tau(x,1)]. \tag{SM2.8}
\end{align*}
\end{theorem}

\begin{theorem}[Theorem 5.2 in \cite{li2021numerical}]
Consider an  operator $\mathcal{L}$ on $\bar{\mathbf{X}}$ in the form of $
    \mathcal{L} f(x) :=  \sum_{i=1}^{J} \zeta_i(x) [f(c_i)-f(c_{i-1})], \, \forall x\in\mathbb{R}, f\in\bar{\mathbf{X}}.$
Then $(\mathcal{I}-\mathcal{L})^{-1}$ exists iff matrix $H$ is non-singular, where $H$ is defined by $H_{ij}:= \delta_{ij}-\zeta_i(c_j)+\zeta_{i+1}(c_j), \, i,j=1,2,...,J$, and $\zeta_{J+1}(x)\equiv0,\forall x\in\mathbb{R}$. 
\end{theorem}

\begin{theorem}[Theorem 21.67 in \cite{hewitt2013real}]
Let $\alpha$ and $\beta$ be any two real-valued  functions on $\mathbb{R}$ of finite variation,  and let $\lambda_\alpha$ and $\lambda_\beta$ be their corresponding Lebesgue-Stieltjes measures. Then $a<b$ in $\mathbb{R}$ implies 
\begin{align*}
    \int_{[a,b]}\frac{\beta(x+)+\beta(x-)}{2}\hat{\mathrm{d}} \lambda_\alpha(x)+\int_{[a,b]}\frac{\alpha(x+)+\alpha(x-)}{2}\hat{\mathrm{d}} \lambda_\beta(x) = \alpha(b+)\beta(b+)-\alpha(a-)\beta(a-).
\end{align*}
\end{theorem}

\begin{theorem}[Theorem 3.4 in \cite{li2021numerical}]
For $\tau \in {T}(\mathbb{R,\mathcal{B}})$ satisfying Condition \ref{A_UBV} (uniformly bounded variation), operator $\mathcal{L}$ in  Definition \ref{Def:K} is compact on $\bar{\mathbf{X}}$ iff Condition \ref{A_MV} (c\`adl\`ag and countable jump discontinuities) holds.
\end{theorem}

\begin{lemma}[Lemma 2.4 in \cite{li2021numerical}] 
    Consider a Banach space $(\mathbf{Y},||\cdot||)$ and a compact linear operator $\mathcal{L}$.  The following statements are equivalent:
    \begin{itemize}
        \item[$(i)$] {\rm Ker}$(\mathcal{I}-\mathcal{L})=\{0\}$,  i.e., $\forall \xi\in\mathbf{Y}$, $f=\xi+\mathcal{L}f$ has at most one solution $f\in \mathbf{Y}$; 
        \item[$(ii)$] {\rm Im}$(\mathcal{I}-\mathcal{L})=\mathbf{Y}$, i.e., $\forall \xi\in\mathbf{Y}$, $f=\xi+\mathcal{L}f$ has at least one solution $f\in \mathbf{Y}$; and
        \item[$(iii)$] $(\mathcal{I}-\mathcal{L})^{-1}$ exists, i.e., $\forall \xi\in\mathbf{Y}$, $f=\xi+\mathcal{L}f$ has a unique solution $f\in \mathbf{Y}$. 
    \end{itemize}
    {A corollary is:} Consider a Banach space $(\mathbf{Y},||\cdot||)$,  linear operators $\mathcal{L}$ and $\{\mathcal{L}^{(r)}\}_{r\in\mathbb{N}}$ on $\mathbf{Y}$. Suppose $\{\mathcal{L}^{(r)}\}_{r\in\mathbb{N}}$ are collectively compact and consistent operators to $\mathcal{L}$. The following statements are equivalent:
    \begin{itemize}
        \item[$(i)$] {\rm Ker}$(\mathcal{I}-\mathcal{L})=\{0\}$,  i.e., $\forall \xi\in\mathbf{Y}$, $f=\xi+\mathcal{L}f$ has at most one solution $f\in \mathbf{Y}$; 
        \item[$(ii)$] {\rm Im}$(\mathcal{I}-\mathcal{L})=\mathbf{Y}$, i.e., $\forall \xi\in\mathbf{Y}$, $f=\xi+\mathcal{L}f$ has at least one solution $f\in \mathbf{Y}$; and
        \item[$(iii)$] $(\mathcal{I}-\mathcal{L})^{-1}$ exists, i.e., $\forall \xi\in\mathbf{Y}$, $f=\xi+\mathcal{L}f$ has a unique solution $f\in \mathbf{Y}$. 
    \end{itemize}
\end{lemma}
\begin{theorem}[Theorem 10.0.1 in \cite{meyn2012markov}]
If the MC is recurrent then it admits a unique (up to constant multiples) invariant measure. The invariant measure is finite (rather than merely $\sigma$-finite) if there exists a petite
set C such that $\sup_{x\in C} \mathbb{E}_x [\tau_C] <\infty$.
\end{theorem}

\begin{theorem}[Theorem 8.3.6 in \cite{meyn2012markov}]
Suppose an MC is $\psi$-irreducible. Then it is recurrent if there exists some petite set $C$ such that with any initial state $x\in C$, the MC returns to $C$ with probability $1$.
\end{theorem}

\begin{comment}
\section{Error bound's continuity with respect to model parameters.}\label{append:parameter continuity}
We assume the model includes a parameter $\mu$ and the approximate MC's transition matrix $Q^{(r)}$ is continuous w.r.t. $\mu$. Then the sensitivity measure $e_2$ in Theorem \ref{Theorem: error bound naive} is also continuous w.r.t. $\mu$. Formally, we have that:
\begin{lemma}
Consider an instance of approximate MC in Approximation \ref{A:Stationary Dist Eq}, if the transition matrix $Q^{(r)}$ is continuous w.r.t. $\mu\in E\subseteq \mathbb{R}$, then $e_2$ is also continuous w.r.t. $\mu$. 

\end{lemma}
\end{comment}

\end{APPENDICES}
%
%   or 
%
% \begin{APPENDICES}
% \section{<Title of Section A>}
% \section{<Title of Section B>}
% etc
% \end{APPENDICES}

% Acknowledgments here
%\section*{Acknowledgments.}
% Enter the text of acknowledgments here

% References here (outcomment the appropriate case) 

% CASE 1: BiBTeX used to constantly update the references 
%   (while the paper is being written).
\bibliographystyle{informs2014} % outcomment this and next line in Case 1
\bibliography{reference} % if more than one, comma separated

% CASE 2: BiBTeX used to generate mypaper.bbl (to be further fine tuned)
%\input{mypaper.bbl} % outcomment this line in Case 2

\end{document}